\documentclass[10pt,a4paper]{article}
\usepackage[latin1]{inputenc}
\usepackage[english]{babel}
\usepackage{amsfonts}
\usepackage{amssymb}
\usepackage{amsmath}
\usepackage{latexsym}
\usepackage[T1]{fontenc}
\usepackage{subfig}
\usepackage{tikz}
\newtheorem{theoreme}{Theorem}[section]
\newtheorem{lemme}[theoreme]{Lemma}
\newtheorem{corollaire}[theoreme]{Corollary}
\newtheorem{proposition}[theoreme]{Proposition}
\newtheorem{remarque}[theoreme]{Remark}

\newtheorem{defi}[theoreme]{Definition}
\newenvironment{preuve}{\emph{Proof} : }{\begin{flushright}$\Box$\end{flushright}}
\newcommand{\F}{\mathbb{F}}

\title{Second weight codewords of generalized Reed-Muller codes}
\author{Elodie Leducq}
\date{} 
\begin{document}
\maketitle
\section{Introduction}

\indent In this paper, we want to characterize the second weight codewords of generalized Reed-Muller codes.
\\\\\indent We first introduce some notations : 
\\\\\indent Let $p$ be a prime number, $n$ a positive integer, $q=p^n$ and $\F_q$ a finite field with $q$ elements.
\\\\\indent If $m$ is a positive integer, we denote by $B_m^q$ the $\F_q$-algebra of the functions from $\F_q^m$ to $\F_q$ and by $\F_q[X_1,\ldots,X_m]$ the $\F_q$-algebra of polynomials in $m$ variables with coefficients in $\F_q$. 

We consider the morphism of $\F_q$-algebras $\varphi: \F_q[X_1,\ldots,X_m]\rightarrow B_m^q$ which associates to $P\in\F_q[X_1,\ldots,X_m]$ the function $f\in B_m^q$ such that $$\textrm{$\forall x=(x_1,\ldots,x_m)\in\F_q^m$, $f(x)=P(x_1,\ldots,x_m)$.}$$ The morphism $\varphi$ is onto and its kernel is the ideal generated by the polynomials $X_1^q-X_1,\ldots,X_m^q-X_m$. So, for each $f\in B_m^q$, there exists a unique polynomial $P\in\F_q[X_1,\ldots,X_m]$ such that the degree of $P$ in each variable is at most $q-1$ and $\varphi(P)=f$. We say that $P$ is the reduced form of $f$ and we define the degree $\deg(f)$ of $f$ as the degree of its reduced form.
 The support of $f$ is the set $\{x\in\F_q^m:f(x)\neq0\}$ and we denote by $|f|$ the cardinal of its support (by identifying canonically $B_m^q$ and $\F_q^{q^m}$, $|f|$ is actually the Hamming weight of $f$).
\\\\\indent For $0\leq r\leq m(q-1)$, the $r$th order generalized Reed-Muller code of length $q^m$ is $$R_q(r,m):=\{f\in B_m^q :\deg(f)\leq r\}.$$
\indent For $1\leq r\leq m(q-1)-2$, the automorphism group of generalized Reed-Muller codes $R_q(r,m)$ is the affine group of $\F_q^m$ (see \cite{charpin_auto}).
\\\\\indent For more results on generalized Reed-Muller codes, we can see for example \cite{delsarte_poids_min}.
\\\\\indent We are now able to give precisely some results about minimum weight codewords and second weight codewords :
\\\\\indent We write $r=t(q-1)+s$, $0\leq t\leq m-1$, $0\leq s\leq q-2$. \\\\\indent In \cite{MR0275989}, interpreting generalized Reed-Muller codes in terms of BCH codes, it is proved that the minimal weight of the generalized Reed-Muller code $R_q(r,m)$ is $(q-s)q^{m-t-1}$. 
\\\\\indent The following theorem gives the minimum weight codewords of generalized Reed-Muller codes and is proved in \cite{delsarte_poids_min} or \cite{Leducq2012581}
\begin{theoreme}\label{poidsmin}Let $r=t(q-1)+s<m(q-1)$, $0\leq s\leq q-2$. The minimal weight codewords of $R_q(r,m)$ are codewords of $R_q(r,m)$ whose support is the union of $(q-s)$ distinct parallel affine subspaces of codimension $t+1$ included in an affine subspace of codimension $t$. \end{theoreme}

In \cite{MR2411119}, Geil proves that the second weight of generalized Reed-Muller codes $R_q((m-1)(q-1)+s,m)$, $1\leq s\leq q-2$ is $q-s+1$ and that the second weight of generalized Reed-Muller codes $R_q(r,m)$, $2\leq r<q$ is $(q-r+1)(q-1)q^{m-2}$. The other cases can be found in the following theorem. Rolland proves all the cases such that $s\neq1$ in \cite{MR2592428}. The case where $s=1$ has been proved by Bruen in \cite{MR2766082} using methods of Erickson (see \cite{erickson1974counting}):

\begin{theoreme}\label{poids2}For $m\geq 3$, $q\geq 3$ and $q\leq r\leq(m-1)(q-1)$ the second weight $W_2$ of the generalized Reed-Muller codes $R_q(r,m)$ satisfies :
\begin{enumerate}\item if $1\leq t\leq m-1$ and $s=0$, $$W_2=2(q-1)q^{m-t-1};$$
\item if $1\leq t \leq m-2$ and $s=1$,
\begin{enumerate}
\item if $q=3$, $W_2=8\times 3^{m-t-2},$
\item if $q\geq4$, $W_2=q^{m-t}$,\end{enumerate}
\item if $1\leq t\leq m-2$ and $2\leq s\leq q-2$, $$W_2=(q-s+1)(q-1)q^{m-t-2}.$$\end{enumerate} \end{theoreme} 

In \cite{MR1384161}, Cherdieu and Rolland prove that the codewords of the second weight of $R_q(s,m)$, $2\leq s\leq q-2$, which are the product of $s$ polynomials of degree 1 are of the following form.

\begin{theoreme}\label{th2.2}Let $m\geq2$, $2\leq s\leq q-2$ and $f\in R_q(s,m)$ such that \\$|f|=(q-s+1)(q-1)q^{m-2}$; we denote by $S$ the support of $f$. Assume that $f$ is the product of $s$ polynomials of degree 1 then either $S$ is the union of $q-s+1$ parallel affine hyperplanes minus their intersection with an affine hyperplane which is not parallel or $S$ is the union of $(q-s+1)$ affine hyperplanes which meet in a common affine subspace of codimension 2 minus this intersection.\end{theoreme}

In \cite{MR2332391}, Sboui proves that the only codewords of $R_q(s,m)$, $2\leq s\leq\frac{q}{2}$ whose weight is $(q-s+1)(q-1)q^{m-2}$ are these codewords.\\

All the results proved in this paper are summarized in Section 2 and their proofs are in the following sections.

\section{Results}
In the following, except when an other affine space is specified, an hyperplane or a subspace is an affine hyperplane or an affine subspace of $\F_q^m$.
\subsection{Case where $t=m-1$ and $s\neq0$}

\begin{theoreme}\label{t=m-1}Let $m\geq 2$, $q\geq 5$, $1\leq s\leq q-4$ and $f\in R_q((m-1)(q-1)+s,m)$ such that $|f|=q-s+1$. Then the support of $f$ is included in a line.\end{theoreme}

\begin{proposition}\label{t=m-1b}Let $m\geq 2$. If $q\geq3$ and $f\in R_q((m-1)(q-1)+q-3,m)$ such that $[f|=4$ or $f\in R_q((m-1)(q-1+q-2,m)$ such that $|f|=3$, then the support of $f$ is included in an affine plane. \end{proposition}

\subsection{Case where $0\leq t\leq m-2$ and $2\leq s \leq q-2$}

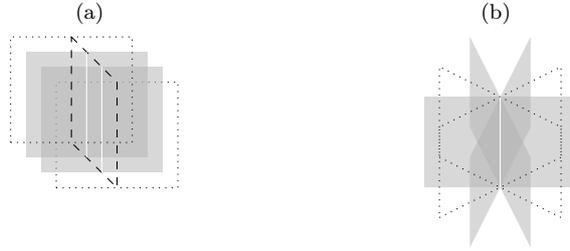
\begin{figure}[!h]
\caption{The possible support for a second weight codeword of $R_4(5,3)$}
\begin{center}\subfloat[]{\label{fig1a}
\begin{tikzpicture}[scale=0.2]
\draw[dotted](3,0)--++(8,0)--++(0,7)--++(-8,0)--cycle;
\fill[color=gray!60, opacity=0.5](2,1)--++(8,0)--++(0,7)--++(-8,0)--cycle;
\fill[color=gray!60, opacity=0.5](1,2)--++(8,0)--++(0,7)--++(-8,0)--cycle;
\draw[dotted](0,3)--++(8,0)--++(0,7)--++(-8,0)--cycle;
\draw[dashed](4,10)--(7,7)--(7,0)--(4,3)--cycle;
\draw[color=white] (5,9)--(5,2);
\draw[color=white] (6,8)--(6,1);
\end{tikzpicture}\label{fig1a}}
\hspace{3cm}
\subfloat[]{\label{fig1b}
\begin{tikzpicture}[scale=0.2]
\fill[color=gray!60, opacity=0.5](0,4)--(10,4)--(10,10)--(0,10)--cycle;
\fill[color=gray!60, opacity=0.5](3,0)--(7,8)--(7,14)--(3,6)--cycle;
\fill[color=gray!60, opacity=0.5](7,0)--(7,6)--(3,14)--(3,8)--cycle;
\draw[dotted](1,2)--(9,6)--(9,12)--(1,8)--cycle;
\draw[dotted](1,6)--(9,2)--(9,8)--(1,12)--cycle;
\draw[color=white](5,4)--(5,10);
\end{tikzpicture}}\end{center}
\label{fig1}
\end{figure}

\begin{theoreme}\label{wpoids2}Let $q\geq4$, $m\geq 2$, $0\leq t\leq m-2$, $2\leq s\leq q-2$. The second weight codewords of $R_q(t(q-1)+s,m)$ are codewords of $R_q(t(q-1)+s,m)$ whose support $S$ is included in an affine subspace of codimension $t$ and either $S$ is the union of $q-s+1$ parallel affine subspaces of codimension $t+1$ minus their intersection with an affine subspace of codimension $t+1$ which is not parallel or $S$ is the union of $(q-s+1)$ affine subspaces of codimension $t+1$ which meet in an affine subspace of codimension $t+2$ minus this intersection (see Figure \ref{fig1}).\end{theoreme}

\subsection{Case where $s=0$}

\begin{theoreme}\label{0} Let $m\geq 2$, $q\geq3$, $1\leq t\leq m-1$. The second weight codewords of $R_q(t(q-1),m)$ are codewords of $R_q(t(q-1),m)$ whose support $S$ is included in an affine subspace of codimension $t-1$ and either $S$ is the union of $2$ parallel affine subspaces of codimension $t$ minus their intersection with an affine subspace of codimension $t$ which is not parallel or $S$ is the union of 2 non parallel affine subspaces of codimension $t$ minus their intersection. \end{theoreme}

\subsection{Case where $0\leq t\leq m-2$ and $s=1$}

\begin{theoreme}\label{s=1}For $q\geq 4$, $m\geq1$, $0\leq t\leq m-1$, if $f\in R_q(t(q-1)+1,m)$ is such that $|f|=q^{m-t}$, the support of $f$ is an affine subspace of codimension $t$. \end{theoreme}

\begin{figure}[!h]
\caption{The support of a second weight codeword of $R_3(3,3)$}
\label{fig2}
\begin{center}\begin{tikzpicture}[scale=0.2]
\fill[color=gray!60, opacity=0.5](0,0)--++(4,2)--++(0,9)--++(-4,-2)--cycle;
\fill[color=gray!60, opacity=0.5](6,0)--++(4,2)--++(0,9)--++(-4,-2)--cycle;
\draw[dotted](3,0)--++(4,2)--++(0,9)--++(-4,-2)--cycle;
\draw[white](0,3)--++(4,3);
\draw[white](6,3)--++(4,3);
\draw[dashed](0,3)--++(6,0);
\draw[dashed](4,6)--++(6,0);
\draw[dashed](0,7)--++(6,0);
\draw[dashed](4,4)--++(6,0);
\draw[white](0,7)--++(4,-3);
\draw[white](6,7)--++(4,-3);
\end{tikzpicture}\end{center}
\end{figure}
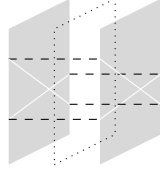

\begin{proposition}\label{1q=3}Let $m\geq3$, $1\leq t\leq m-2$ and $f\in R_3(2t+1,m)$ such that $|f|=8.3^{m-t-2}$. We denote by $S$ the support of $f$. Then $S$ is included in $A$ an affine subspace of dimension $m-t+1$, $S$ is the union of two parallel hyperplanes of $A$ minus their intersection with two non parallel hyperplanes of $A$ (see Figure \ref{fig2}).\end{proposition}

\section{Some tools}

The following lemma and its corollary are proved in \cite{delsarte_poids_min}.

\begin{lemme}\label{DGMW1}Let $m\geq1$, $q\geq2$, $f\in B_m^q$ and $a\in\F_q$. If for all $(x_2,\ldots,x_m)$ in $\F_q^{m-1}$, $f(a,x_2,\ldots,x_m)=0$ then for all $(x_1,\ldots,x_m)\in\F_q^m$, $$f(x_1,\ldots,x_m)=(x_1-a)g(x_1,\ldots,x_m)$$ with $\deg_{x_1}(g)\leq\deg_{x_1}(f)-1$.\end{lemme}

\begin{corollaire}Let $m\geq1$, $q\geq2$, $f\in B_m^q$ and $a\in\F_q$. If for all $(x_1,\ldots,x_m)$ in $\F_q^m$ such that $x_1\neq a$, $f(x_1,\ldots,x_m)=0$ then for all $(x_1,\ldots,x_m)\in\F_q^m$, $f(x_1,\ldots,x_m)=(1-(x_1-a)^{q-1})g(x_2,\ldots,x_m)$.\end{corollaire} 

\begin{lemme}\label{inter}Let $q\geq3$, $m\geq3$, and $S$ be a set of points of $\F_q^m$ such that $\#S=u.q^{n}<q^m$, with $u\not\equiv0\mod q$. Assume that for all hyperplanes $H$ either $\#(S\cap H)=0$ or $\#(S\cap H)=v.q^{n-1}$, $v<u$ or $\#(S\cap H)\geq u.q^{n-1}$ Then there exists $H$ an affine hyperplane such that $S$ does not meet $H$ or such that $\#(S\cap H)=vq^{n-1}$.\end{lemme}

\begin{preuve}Assume that for all $H$ hyperplane, $S\cap H\neq \emptyset$ and $\#(S\cap H)\neq vq^{n-1}$. Consider an affine hyperplane $H$; then for all $H'$ hyperplane parallel to $H$, $\#(S\cap H')\geq u.q^{n-1}.$ 
Since $u.q^{n}=\#S=\displaystyle\sum_{H'//H}\#(S\cap H')$, we get that for all  $H$ hyperplane, $\#(S\cap H)=u.q^{n-1}$.
\\Now  consider $A$ an affine subspace of codimension 2 and the $(q+1)$ hyperplanes through $A$. These hyperplanes intersect only in $A$ and their union is equal to $\F_q^m$. So $$uq^{n}=\#S=(q+1)u.q^{n-1}-q\#(S\cap A).$$
Finally we get a contradiction if $n=1$. Otherwise, $\#(S\cap A)=u.q^{n-2}$. Iterating this argument, we get that for all $A$ affine subspace of codimension $k\leq n$, $\#(S\cap A)=u.q^{n-k}$. 
\\Let $A$ be an affine subspace of codimension $n+1$ and  $A'$ an affine subspace of codimension $n-1$ containing $A$. We consider the $(q+1)$ affine subspace of codimension $n$ containing $A$ and  included in $A'$, then $$u.q=\#(S\cap A')=(q+1)u-q\#(S\cap A)$$ which is absurd since $\#(S\cap A)$ is an integer and $u\not\equiv 0\mod q$. So there exists $H_0$ an hyperplane such that $\#(S\cap H_0)=vq^{n-1}$ or $S$ does not meet $H_0$. \end{preuve}

\section{Case where $t=m-1$ and $s\neq0$}

\subsection{Proof of Theorem \ref{t=m-1}}

Let $\omega_1$, $\omega_2\in S$ and $H$ an affine hyperplane containing $\omega_1$ and $\omega_2$. Assume $S\cap H\neq S$. We have $\#S=q-s+1\leq q$ and $\omega_1$, $\omega_2\in S\cap H$, so there exists an affine hyperplane parallel to $H$ which does not meet $S$. By applying an affine transformation, we can assume that $x_1=0$ is an equation of $H$ and we denote by $H_a$ the affine hyperplane parallel to $H$ of equation $x_1=a$, $a\in\F_q$. Let $I:=\{a\in\F_q:S\cap H_a=\emptyset\}$ and denote by $k:=\# I$; $s\leq k\leq q-2$. Let $c\not\in I$, we define 
$$\forall x=(x_1,\ldots,x_m)\in\F_q^m, \ f_c(x)=f(x)\prod_{a\not\in I,a\neq c}(x_1-a)$$
that is to say $f_c$ is a function in $B_m^q$ such that its support is $S\cap H_c$. Since $c\not\in I$, $f_c$ is not identically zero. Then $|f|=\displaystyle\sum_{c\not\in I}|f_c|$ and we consider two cases.
\begin{itemize}\item Assume $k>s$.
\\Then the reduced form of $f_c$ has degree at most $(m-1)(q-1)+q-1+s-k$ and $|f_c|\geq k-s+1$. Then, $$(q-s+1)=|f|=\sum_{c\not\in I}|f_c|\geq(q-k)(k-s+1)$$ which gives $$1\geq(q-1-k)(k-s)$$ this is possible if and only if $k=q-2=s+1$ and we get a contradiction since $s\leq q-4$.
\item Assume that $k=s$.
\\Then $S$ meets $(q-s-1)$ affine hyperplanes parallel to $H$ in 1 point and $H$ in 2 points. Consider the function $g$ in $B_m^q$ defined by $$\forall x=(x_1,\ldots,x_m)\in\F_q^m, \ g(x)=x_1f(x).$$ The reduced form of $g$ has degree at most $(m-1)(q-1)+s+1$ and $$|g|=(q-s-1).$$ So $g$ is a minimum weight codeword of $R_q((m-1)(q-1)+s+1,m)$ and its support is included in a line. This line is not included in $H$. So consider $H_1$ an affine hyperplane which contains this line but does not contain both $\omega_1$ and $\omega_2$. Then $S\cap H_1\neq S$ and $H_1$ contains at least 3 points of $S$ since $s\leq q-4$ which gives a contradiction by applying the previous argument to $H_1$.\end{itemize}
So $S$ is included in all affine hyperplanes through $\omega_1$ and $\omega_2$ which gives the result.

\subsection{Proof of Theorem \ref{t=m-1b}}

\begin{itemize}
\item If $f\in R_q((m-1)(q-1)+q-2,m)$ is such that $|f|=3$, we have the result since 3 points are always included in an affine plane.
\item Assume $f\in R_q((m-1)(q-1)+q-3,m)$ is such that $|f|=4$. 
\\Let $a$, $b$, $c$, $d\in\F_q^*$ and $\omega^{(a)}=(\omega_1^{(a)},\ldots,\omega_m^{(a)})$, $\omega^{(b)}=(\omega_1^{(b)},\ldots,\omega_m^{(b)})$, $\omega^{(c)}=(\omega_1^{(c)},\ldots,\omega_m^{(c)})$, $\omega^{(d)}=(\omega_1^{(d)},\ldots,\omega_m^{(d)})$ 4 distinct points of $\F_q^m$ such that $\forall x=(x_1,\ldots,x_m)\in\F_q^m,$ \begin{align*}f(x)=&a\prod_{i=1}^m\left(1-(x_i-\omega_i^{(a)})^{q-1}\right)+b\prod_{i=1}^m\left(1-(x_i-\omega_i^{(b)})^{q-1}\right) \\&+c\prod_{i=1}^m\left(1-(x_i-\omega_i^{(c)})^{q-1}\right)+d\prod_{i=1}^m\left(1-(x_i-\omega_i^{(d)})^{q-1}\right).\end{align*}
So, \begin{align*}f(x)=&(-1)^m(a+b+c+d)\prod_{i=1}^mx_i^{q-1} \\&+(-1)^{m-1}\sum_{i=1}^m(a\omega_i^{(a)}+b\omega_i^{(b)}+c\omega_i^{(c)}+d\omega_i^{(d)})x_i^{q-2}\prod_{j\neq i}x_j^{q-1}+r\end{align*}
with $\deg(r)\leq(m-1)(q-1)+q-3$.
Since $f\in R_q((m-1)(q-1)+q-3,m)$, $$\left\{\begin{array}{l}a+b+c+d=0\\a\omega^{(a)}+b\omega^{(b)}+c\omega^{(c)}+d\omega^{(d)}=0\end{array}\right..$$
So, $a\overrightarrow{\omega^{(d)}\omega^{(a)}}+b\overrightarrow{\omega^{(d)}\omega^{(b)}}+c\overrightarrow{\omega^{(d)}\omega^{(c)}}=\overrightarrow{0}$ which gives the result.
\end{itemize}

\begin{remarque}In both cases we cannot prove that the support of $f$ is included in a line. Indeed, 
\begin{itemize}\item Let $\omega_1$, $\omega_2$, $\omega_3$ 3 points of $\F_q^m$ not included in a line. For $q\geq 3$ we can find $a$, $b\in\F_q^*$ such that $a+b\neq0$. Let $ f=a\mathrm{1}_{\omega_1}+b\mathrm{1}_{\omega_2}-(a+b)\mathrm{1}_{\omega_3}$ where  for $\omega\in\F_q^m$, $\mathrm{1}_{\omega}$ is the function from $\F_q^m$ to $\F_q$ such that $\mathrm{1}_{\omega}(\omega)=1$ and $\mathrm{1}_{\omega}(x)=0$ for all $x\neq \omega$. Then, since $\displaystyle\sum_{x\in\F_q^m}f(x)=a+b-(a+b)=0$, $f\in R_q((m-1)(q-1)+q-2,m)$.
\item Let $\omega_1$, $\omega_2$, $\omega_3$ 3 points of $\F_q^m$ not included in a line and set $$\omega_4=\omega_1+\omega_2-\omega_3.$$ Then $f=\mathrm{1}_{\omega_1}+\mathrm{1}_{\omega_2}-\mathrm{1}_{\omega_3}-\mathrm{1}_{\omega_4}\in R_q((m-1)(q-1)+q-3,m)$.\end{itemize}\end{remarque}

\section{Case where $0\leq t\leq m-2$ and $2\leq s \leq q-2$}

\subsection{Case where $t=0$}

In this subsection, we write $r=a(q-1)+b$ with $0\leq a\leq m-1$ and $0<b\leq q-1$. 

\begin{lemme}\label{cle}Let $q\geq3$, $m\geq2$, $0\leq a\leq m-2$, $2\leq b\leq q-1$ and \\$f\in R_q(a(q-1)+b,m)$ such that $|f|=(q-b+1)(q-1)q^{m-a-2}$; we denote by $S$ the support of $f$. If $H$ is an affine hyperplane of $\F_q^m$ such that $S\cap H\neq\emptyset$ and $S\cap H\neq S$ then either $S$ meets all affine hyperplanes parallel to $H$ or $S$ meets $q-b+1$ affine hyperplanes parallel to $H$ in $(q-1)q^{m-a-2}$ points or $S$ meets $q-1$ affine hyperplanes parallel to $H$ in $(q-b+1)q^{m-a-2}$ points. \end{lemme}

\begin{preuve}By applying an affine transformation, we can assume that $x_1=0$ is an equation of $H$ and consider the $q$ affine hyperplanes $H_w$ of equation $x_1=w$, $w\in\F_q$, parallel to $H$. Let $I:=\{w\in\F_q:S\cap H_w=\emptyset\}$ and denote by $k:=\# I$. Assume that $k\geq1$. Since $S\cap H\neq \emptyset$ and $S\cap H\neq S$, $k\leq q-2$. For all $c\in\F_q$, $c\not\in I$, we define $$\forall x=(x_1,\ldots,x_n)\in\F_q^m, \ f_c(x)=f(x)\prod_{w\in\F_q,w\neq c,w\not\in I}(x_1-w).$$
 
\begin{itemize}\item Assume $b<k$.
\\Then $2\leq q-1+b-k\leq q-2$ and for all $c\not\in I$, the reduced form of $f_c$ has degree at most $a(q-1)+q-1+b-k$. So $|f_c|\geq (k-b+1)q^{m-a-1}$. Hence $$(q-1)(q-b+1)q^{m-a-2}\geq(q-k)(k-b+1)q^{m-a-1}$$ which means that $(b-k)q(q-k-1)+b-1\geq0$. However $(b-k)\leq -1$ and $q-k-1\geq 1$ so $(b-k)q(q-k-1)+b-1<0$ which gives a contradiction. 
\item Assume $b\geq k$.
\\Then $0\leq b-k\leq q-2$ and for all $c\not\in I$, the reduced form of $f_c$ has degree at most $(a+1)(q-1)+b-k$. So $|f_c|\geq (q-b+k)q^{m-a-2}$. Hence $$(q-1)(q-b+1)q^{m-a-2}\geq (q-k)(q-b+k)q^{m-a-2}$$ with equality if and only if for all $c\not\in I$, $|f_c|=(q-b+k)q^{m-a-2}$. Finally, we obtain that $(k-1)(k-b+1)\geq0$ which is possible if and only if  $k=1$ or $1\geq b-k\geq0$. Now, we have to show that $k=s$ is impossible to prove the lemma. If $b=q-1$, since $k\leq q-2$, we have the result. Assume that $b\leq q-2$ and $b=k$. Then, for all $c\not\in I$, $f_c\in R_q((a+1)(q-1),m)$. The minimum weight of $R_q((a+1)(q-1),m)$ is $q^{m-a-1}$ and its second weight is $2(q-1)q^{m-a-2}$. We denote by $N_1:=\#\{c\not\in I:|f_c|=q^{m-a-1}\}$. Since $k=b$, $N_1\leq q-b$. Furthermore, we have $$(q-b+1)(q-1)q^{m-a-2}\geq N_1q^{m-a-1}+(q-b-N_1)2(q-1)q^{m-a-2}$$ which means that $N_1\geq \frac{(q-1)(q-b-1)}{q-2}>q-b-1$. Finally, $N_1=q-b$ and for all $c\not\in I$, $|f_c|=q^{m-a-1}$. However $(q-1)(q-b+1)q^{m-a-2}>(q-b)q^{m-a-1}$ which gives a contradiction.\end{itemize}\end{preuve}

\begin{lemme}\label{m=2}For $m=2$, $q\geq3$, $2\leq b\leq q-1$. The second weight codewords of $R_q(b,2)$ are codewords of $R_q(b,2)$ whose support $S$ is the union of $q-b+1$ parallel lines minus their intersection with a line which is not parallel or $S$ is the union of $(q-b+1)$ lines which meet in a point minus this point.\end{lemme}

\begin{preuve}To prove this lemma, we use some results on blocking sets proved by Erickson in \cite{erickson1974counting} and Bruen in \cite{MR2766082}. All these results are recalled in the Appendix of this paper. By Theorem \ref{th2.2}, which is also true for $b=q-1$ (see \cite[Lemma 3.12]{erickson1974counting}), it is sufficient to prove that $f\in R_q(b,2)$ such that $|f|=(q-b+1)(q-1)$  is the product of linear factors. 

Let $f\in R_q(b,2)$ such that $|f|\leq (q-b+1)(q-1)=q(q-b)+b-1.$ We denote by $S$ its support. Then, $S$ is not a blocking set of order $(q-b)$ of $\F_q^2$ (Theorem \ref{4.14}) and $f$ has a linear factor (Lemma \ref{4.2}). 

We proceed by induction on $b$. If $b=2$ and $f\in R_q(b,2)$ is such that $|f|\leq (q-b+1)(q-1)$, then $f$ has a linear factor and by Lemma \ref{DGMW1} $f$ is the product of 2 linear factors. Assume that if $f\in R_q(b-1,2)$ is such that $|f|\leq (q-b+2)(q-1)$ then $f$ is a product of linear factors. Let $f\in R_q(b,2)$ such that $|f|\leq (q-b+1)(q-1)$; then $f$ has a linear factor. By applying an affine transformation, we can assume that for all $(x,y)\in\F_q^2$, $f(x,y)=y\widetilde{f}(x,y)$ with $\deg(\widetilde{f})\leq b-1$.  So, $L$ the line of equation $y=0$ does not meet $S$ the support of $f$. Since $(q-b+1)(q-1)>q$, $S$ is not included in a line and by Lemma \ref{cle}, either $S$ meets $(q-b+1)$ lines parallel to $L$ in $(q-1)$ points or $S$ meets $(q-1)$ lines parallel to $L$ in $(q-b+1)$ points. 

In the first case, by Lemma \ref{DGMW1}, we can write for all $(x,y)\in\F_q^2$, $$f(x,y)=y(y-a_1)\ldots(y-a_{b-2})g(x,y)$$ where $a_i$, $1\leq i\leq q-2$ are $q-2$ distinct elements of $\F_q^*$ and $\deg(g)\leq 1$ which gives the result.

In the second case, we denote by $a\in\F_q$ the coefficient of $x^{s-1}$ in $\widetilde{f}$. Then for any $\lambda\in\F_q^*$, since $S$ meets all lines parallel to $L$ but $L$ in $q-s+1$ points, we get for all $x\in\F_q$, $$f(x,\lambda)=a\lambda(x-a_1(\lambda))\ldots(x-a_{b-1}(\lambda))$$
So there exists $a_1,\ldots a_{b-1}\in\F_q[Y]$ of degree at most $q-1$ such that for all $(x,y)\in\F_q^2$,
$$f(x,y)=ay(x-a_1(y))\ldots(x-a_{b-1}(y)).$$ Then for all $x\in\F_q$, $$\widetilde{f}_0(x)=\widetilde{f}(x,0)=a(x-a_1(0))\ldots(x-a_{b-1}(0))$$ and $|\widetilde{f}_0|\leq q-1$. So, $$|\widetilde{f}|\leq |f|+|\widetilde{f}_0|\leq (q-b+2)(q-1).$$
By recursion hypothesis, $\widetilde{f}$ is the product of linear factors which finishes the proof of Lemma \ref{m=2}.\end{preuve}

\begin{proposition}\label{t=0}For $m\geq2$, $q\geq3$, $2\leq b\leq q-1$. The second weight codewords of $R_q(b,m)$ are codewords of $R_q(b,m)$ whose support $S$ is the union of $q-b+1$ parallel hyperplanes minus their intersection with an affine hyperplane which is not parallel or $S$ is the union of $(q-b+1)$ hyperplanes which meet in an affine subspace of codimension 2 minus this intersection.\end{proposition}

\begin{preuve} We say that we are in configuration $A$ if $S$ is the union of $q-b+1$ parallel hyperplanes minus their intersection with an affine hyperplane which is not parallel (see Figure \ref{fig1a}) and that we are in configuration $B$ if $S$ is the union of $(q-b+1)$ hyperplanes which meet in an affine subspace of codimension 2 minus this intersection (see Figure \ref{fig1b}).

We prove this proposition by induction on $m$. The Lemma \ref{m=2} proves the case where $m=2$. Assume that $m\geq3$ and that second weight codeword of $R_q(b,m-1)$, $2\leq b\leq q-1$ are of type $A$ or type $B$. Let $f\in R_q(b,m)$ such that $|f|=(q-1)(q-b+1)q^{m-2}$ and we denote by $S$ its support.

\begin{itemize}\item Assume that $S$ meets all affine hyperplanes. \\

Then, by Lemma \ref{inter}, there exists an affine hyperplane $H$ such that \\$\#(S\cap H)=(q-b)q^{m-2}$. By applying an affine transformation, we can assume that $x_1=0$ is an equation of $H$. We denote by $\mathrm{1}_{H}$ the function in $B_m^q$ such that $$\forall x=(x_1,\ldots,x_m)\in\F_q^m,\ \mathrm{1}_{H}(x)=1-x_1^{q-1}$$ then the reduced form $f.\mathrm{1}_{H}$ has degree at most $(t+1)(q-1)+s$ and the support of $f.\mathrm{1}_{H}$ is $S\cap H$ so $S\cap H$ is the support of a minimal weight codeword of $R_q(q-1+b,m)$ and $S\cap H$ is the union of $(q-b)$ parallel affine subspaces of codimension 2. Consider $P$ an affine subspace of codimension 2 included in $H$ such that $\#(S\cap P)=(q-b)q^{m-3}$. Assume that there are at least 2 hyperplanes through $P$ which meet $S$ in $(q-b)q^{m-2}$ points. Then, there exists $H_1$ an affine hyperplane through $P$ different from $H$ such that $\#(S\cap H_1)=(q-b)q^{m-2}$. So, $S\cap H_1$ is the union of $(q-b)$ parallel affine subspaces of codimension 2. Consider $G$ an affine hyperplane which contains $Q$ an affine subspace of codimension 2 included in $H$ which does not meet $S$ and the affine subspace of codimension 2 included in $H_1$ which meets $Q$ but not $S$ (see Figure \ref{fig3}).
\begin{figure}[!h]
\caption{}
\label{fig3}
\begin{center}\begin{tikzpicture}[scale=0.2]
\draw (0,0)--(5,2)--(5,11)--(0,9)--cycle;
\draw (1,2/5)--(1,9+2/5);
\draw (2,4/5)--(2,9+4/5);
\draw[dotted] (3,6/5)--(3,9+6/5);
\draw[dotted] (4,8/5)--(4,9+8/5);
\draw (0,9) node[above left]{$H$};
\draw[dashed](0,4)--++(5,2);
\draw (0,4) node[left]{$P$};
\draw (0,4)--(12,4)--++(5,2)--++(-12,0);
\draw (17,6) node[above right]{$H_1$};
\draw (1,4+2/5)--++(12,0);
\draw (2, 4+4/5)--++(12,0);
\draw[dotted] (3,4+6/5)--++(12,0);
\draw[dotted] (4,4+8/5)--++(12,0);
\draw (4,8/5)--++(12,0)--++(0,9)--++(-12,0);
\draw (16,9+8/5) node[above right]{$G$};
\draw (4,9+8/5) node[above]{$Q$};
\end{tikzpicture}\end{center}
\end{figure}
 
By applying an affine transformation, we can assume that $x_m=\lambda$, $\lambda\in\F_q$ is an equation of an hyperplane parallel to $G$. For all $\lambda\in\F_q$, we define $f_{\lambda}\in B_{m-1}^q$ by $$\forall(x_1,\ldots,x_{m-1})\in\F_q^{m-1},\qquad f_{\lambda}(x_1,\ldots,x_{m-1})=f(x_1,\ldots,x_{m-1},\lambda).$$
If all hyperplanes parallel to $G$ meets $S$ in $(q-b+1)(q-1)q^{m-3}$ then for all $\lambda\in\F_q$, $f_{\lambda}$ is a second weight codeword of $R_q(b,m-1)$ and its support is of type $A$ or $B$. We get a contradiction if we consider an hyperplane parallel to $G$ which meets $S\cap H$ and $S\cap H_1$. So, there exits $G_1$ an hyperplane parallel to $G$ which meets $S$ in $(q-b)q^{m-2}$ points and $S\cap G_1$ is the union of $(q-b)$ parallel affine subspaces of codimension 2 which is a contradiction. Then for all $H'$ hyperplane through $P$ different from $H$ $\#(S\cap H')\geq(q-1)(q-b+1)q^{m-3}$. Furthermore, $$(q-b)q^{m-2}+q.(q-1)(q-b+1)q^{m-3}-q.(q-b)q^{m-3}=(q-1)(q-b+1)q^{m-2}.$$ Finally, by applying the same argument to all affine hyperplanes of codimension 2 included in $H$ parallel to $P$, we get $q$ parallel hyperplanes $(G_i)_{1\leq i\leq q}$ such that for all $1\leq i\leq q$, $\#(S\cap G_i)=(q-b+1)(q-1)q^{m-3}$ and $\#(S\cap G_i\cap H)=(q-b)q^{m-3}$. Then by recursion hypothesis, $S\cap G_i$ is either of type $A$ or of type $B$. 

If there exists $i_0$ such that $S\cap G_{i_0}$ is of type $A$. Consider $F$ an affine hyperplane containing $R$ an affine subspace of codimension 2 included in $H$ which does not meet $S$ and the affine subspace of codimension 2 included in $G_{i_0}$ which does not meets $S$ but meets $R$. If for all $F'$ hyperplane parallel to $F$, $\#(S\cap F')>(q-b)q^{m-2}$ then $\#(S\cap F')= (q-1)(q-b+1)q^{m-3}$. So $S\cap F'$ is the support of a second weight codeword of $R_q(b,m-1)$ and is either of type $A$ or of type $B$ which is absurd is we consider an hyperplane parallel to $F$ which  meets $S\cap H$. So there exits $F_1$ an affine hyperplane parallel to $F$ which meets $S$ in $(q-b)q^{m-2}$ points. So $S\cap F_1$ is the union of $(q-s)$ parallel affine subspaces of codimension 2 which is absurd since $S\cap G_{i_0}$ is of type $A$ (see Figure \ref{fig4}). 
\begin{figure}[!h]
\caption{}
\label{fig4}
\begin{center}\begin{tikzpicture}[scale=0.2]
\draw (0,0)--(5,2)--(5,11)--(0,9)--cycle;
\draw (1,2/5)--(1,9+2/5);
\draw (2,4/5)--(2,9+4/5);
\draw[dotted] (3,6/5)--(3,9+6/5);
\draw[dotted] (4,8/5)--(4,9+8/5);
\draw (0,9) node[above left]{$H$};
\draw[dashed](0,4)--++(5,2);
\draw (0,4)--(12,4)--++(5,2)--++(-12,0);
\draw (17,6) node[above right]{$G_{i_0}$};
\draw (1,4+2/5)--++(11+5/2,1-2/5);
\draw (2,4+4/5)--++(11+5/2,1-2/5);
\draw (3,4+6/5)--++(11+5/2,1-2/5);
\draw[dotted] (4,4+8/5)--++(9,2/5);
\draw[dashed] (3,4+6/5)--(9,4);
\draw (3,9+6/5)--++(6,-6/5)--++(0,-9)--++(-6,6/5);
\draw (9,9) node[above right]{$F$};
\draw (3,9+6/5) node[above]{$R$};
\end{tikzpicture}\end{center}
\end{figure}

If for all $1\leq i\leq q$, $S\cap G_i$ is of type $B$. Let $H_1$ be the affine hyperplane parallel to $H$ which contains the affine subspace of codimension 3 intersection of the affine subspaces of codimension 2 of $S\cap G_1$. We consider $R$ an affine subspace of codimension 2 included in $H$ which does not meet $S$. Then there is $(q-b+1)$ affine hyperplanes through $R$ which meet $S\cap G_1$ in $(q-b)q^{m-3}$. However, if we denote by $k$ the number of hyperplanes through $R$ which meet $S$ in $(q-b)q^{m-2}$ points, we have $$k(q-b)q^{m-2}+(q+1-k)(q-1)(q-b+1)q^{m-3}\leq (q-1)(q-b+1)q^{m-2}$$ which implies that $k\geq q-b+2$. For all $H'$ hyperplane through $R$ such that $\#(S\cap H')=(q-b)q^{m-2}$, $S\cap H'$ is the union of $(q-b)$ affine subspaces of codimension 2 parallel to $R$ and then $\#(S\cap H'\cap G_1)=(q-b)q^{m-3}$ which is absurd (see Figure \ref{fig5}). 
\begin{figure}[!h]
\caption{}
\label{fig5}
\begin{center}\begin{tikzpicture}[scale=0.2]
\draw (0,0)--(5,2)--(5,11)--(0,9)--cycle;
\draw (1,2/5)--(1,9+2/5);
\draw (2,4/5)--(2,9+4/5);
\draw[dotted] (3,6/5)--(3,9+6/5);
\draw[dotted] (4,8/5)--(4,9+8/5);
\draw (0,9) node[above left]{$H$};
\draw (3,6/5+9) node[above]{$R$};
\draw[dashed](0,4)--++(5,2);
\draw (0,4)--(12,4)--++(5,2)--++(-12,0);
\draw (17,6) node[above right]{$G_{1}$};
\draw (2,4+4/5)--++(12,0);
\draw (1,4+2/5)--(15,4+6/5);
\draw[dotted] (3,4+6/5)--(13,4+2/5);
\draw[dotted] (4,4+8/5)--(12,4);
\draw (6,4)--(11,6);
\draw[color=white] (8,4+4/5) node {$\bullet$} ;
\draw (6,0)--++(5,2)--++(0,9)--++(-5,-2)--cycle;
\draw (11,11) node[above right]{$H_1$};
\end{tikzpicture}\end{center}
\end{figure}

\item So, there exists $H$ an affine hyperplane such that $H$ does not meet $S$.\\

Then, by Lemma \ref{cle}, either $S$ meets $(q-1)$ hyperplanes parallel to $H$ in $(q-b+1)q^{m-2}$ points or $S$ meets $(q-b+1)$ hyperplanes parallel to $H$ in $(q-1)q^{m-2}$ points.

If $S$ meets $(q-b+1)$ hyperplanes parallel to $H$ in $(q-1)q^{m-2}$ points, then , for all $H'$ hyperplane parallel to $H$ such that $S\cap H'\neq \emptyset$, $S\cap H'$ is the support of a minimal weight codeword of $R_q(q,m)$ and is the union of $(q-1)$ parallel affine subspaces of codimension 2. Let $H'$ be an affine hyperplane parallel to $H$ such that $S\cap H'\neq \emptyset$. We denote by $P$ the affine subspace of codimension 2 of $H'$ which does not meet $S$. Consider $H_1$ an affine hyperplane which contains $P$ and a point not in $S$ of an affine hyperplane $H"$ parallel to $H$ which meets $S$. Then $$\#(H_1\setminus S)\geq bq^{m-2}+1.$$ However, if $S\cap H_1\neq \emptyset$, $\#(H_1\setminus S)\leq bq^{m-2}$. So, $S\cap H_1=\emptyset$ and we are in configuration $A$.

If $S$ meets $(q-1)$ hyperplanes parallel to $H$ in $(q-b+1)q^{m-2}$ points. Then for all $H'$ parallel to $H$ different from $H$, $S\cap H'$ is the support of a minimal weight codeword of $R_q((q-1)+b-1,m)$ and is the union of $(q-b+1)$ parallel affine subspaces of codimension 2. Let $H_1$ be an affine hyperplane parallel to $H$ different from $H$ and consider $P$ an affine subspace of codimension 2 included in $H_1$ such that $$\#(S\cap P)=(q-b+1)q^{m-3}.$$ Assume that there exists $H_2$ an affine hyperplane through $P$ such that $\#(S\cap H_2)=(q-b)q^{m-2}$. Then $S\cap H_2$ is the support of a minimal weight  codeword of $R_q(q-1+b,m)$ and is the union of $(q-b)$ parallel affine subspaces of codimension 2 which is absurd since $S\cap H_2$ meets $H_1$ in $S\cap P$ (see Figure \ref{fig6}).\\

\begin{figure}[!h]
\caption{}
\label{fig6}
\begin{center}\begin{tikzpicture}[scale=0.2]
\draw (0,0)--(5,2)--(5,11)--(0,9)--cycle;
\draw (1,2/5)--(1,9+2/5);
\draw (2,4/5)--(2,9+4/5);
\draw (3,6/5)--(3,9+6/5);
\draw[dotted] (4,8/5)--(4,9+8/5);
\draw[dashed](0,4)--++(5,2);
\draw (0,4) node[left]{$P$};
\draw (0,9) node[above left]{$H_1$};
\draw (0,4)--(12,4)--++(5,2)--++(-12,0);
\draw (17,6) node[above right]{$H_2$};
\draw (1,4+2/5)--++(12,0);
\draw (2,4+4/5)--++(12,0);
\draw[dotted] (3,4+6/5)--++(12,0);
\draw[dotted] (4,4+8/5)--++(12,0);
\end{tikzpicture}\end{center}
\end{figure}

Then, for all $H'$ through $P$ $\#(S\cap H')\geq (q-1)(q-b+1)q^{m-3}$. Furthermore, $$(q-b+1)q^{m-2}+q.(q-1)(q-b+1)q^{m-3}-q.(q-b+1)q^{m-3}=(q-1)(q-b+1)q^{m-2}.$$ So for all $H'$ hyperplane through $P$ different from $H_1$, $$\#(S\cap H')=(q-1)(q-b+1)q^{m-3}.$$ By applying the same argument to all affine subspaces of codimension 2 included in $H_1$ parallel to $P$, we get $q$ parallel hyperplanes $(G_i)_{1\leq i\leq q}$ such that for all $1\leq i \leq q$, $\#(S\cap G_i)=(q-b+1)(q-1)q^{m-3}$ and $\#(S\cap G_i\cap H_1)=(q-s+1)q^{m-3}$. By recursion hypothesis, for all $1\leq i\leq q$, either $S\cap G_i$ is of type $A$ or $S\cap G_i$ is of type $B$. 

Assume that there exists $i_0$ such that $S\cap G_{i_{0}}$ is of type $A$. Consider $F$ an affine hyperplane containing $Q$ an affine subspace of codimension 2 included in $H_1$ which does not meet $S$ and the affine subspace of codimension 2 included in $G_{i_0}$ which does not meets $S$ but meets $Q$. Assume that $S$ meets all hyperplanes parallel to $F$ in at least $(q-b)q^{m-t-2}$. If for all $F'$ parallel to $F$, $\#(S\cap F')>(q-b)q^{m-2}$ then $$\#(S\cap F')\geq (q-1)(q-b+1)q^{m-3}.$$ So $S\cap F'$ is the support of a second weight codeword of $R_q(b,m-1)$ and is either of type $A$ or of type $B$ which is absurd is we consider an hyperplane parallel to $F$ which  meets $S\cap H_1$ and $S\cap G_{i_0}$. So, there exits $F_1$ an affine hyperplane parallel to $F$ such that $\#(S\cap F_1)=(q-b)q^{m-2}$. Then, $S\cap F_1$ is the union of $(q-b)$ parallel affine subspaces of codimension 2, which is absurd. 
Finally, there exists an affine hyperplane parallel to $F$ which does not meet $S$. By Lemma \ref{cle}, either $S$ meets $(q-b+1)$ hyperplanes parallel to $F$ in $(q-1)q^{m-2}$ points and we have already seen that in this case $S$ is of type $A$ or $S$ meets $(q-1)$ hyperplanes parallel to $F$ in $(q-b+1)q^{m-2}$ points. In this case, for all $F'$ parallel to $F$ such that $S\cap F'\neq \emptyset$, $S\cap F'$ is the support of a minimal weight codeword of $R_q(q-1+b-1,m)$ and is the union of $q-b+1$ parallel affine subspaces of codimension 2, which is absurd since $S\cap G_{i_0}$ is of type $A$ (see Figure \ref{fig7}).

\begin{figure}[!h]
\caption{}
\label{fig7}
\begin{center}\begin{tikzpicture}[scale=0.2]
\draw (0,0)--(5,2)--(5,11)--(0,9)--cycle;
\draw (1,2/5)--(1,9+2/5);
\draw (2,4/5)--(2,9+4/5);
\draw[dotted] (3,6/5)--(3,9+6/5);
\draw[dotted] (4,8/5)--(4,9+8/5);
\draw (0,9) node[above left]{$H_1$};
\draw[dashed](0,4)--++(5,2);
\draw (0,4) node[left]{$P$};
\draw (0,4)--(12,4)--++(5,2)--++(-12,0);
\draw (17,6) node[above right]{$G_{i_0}$};
\draw (1,4+2/5)--++(12,0);
\draw (2, 4+4/5)--++(12,0);
\draw[dotted] (3,4+6/5)--++(12,0);
\draw[dotted] (4,4+8/5)--++(12,0);
\draw (4,8/5)--++(12,0)--++(0,9)--++(-12,0);
\draw (16,9+8/5) node[above right]{$F$};
\draw[dashed] (4,4+8/5)--(11,4);
\end{tikzpicture}\end{center}
\end{figure}

Now, assume that for all $1\leq i\leq q$, $G_i\cap S$ is of type $B$. Let $Q$ be an affine subspace of codimension 2 included in $H_1$ which does not meets $S$. Assume that $S$ meets all affine hyperplanes through $Q$ and denote by $k$ the number of these hyperplanes which meet $S$ in $(q-b)q^{m-2}$ points. Then, $$k(q-b)q^{m-2}+(q+1-k)(q-1)(q-b+1)q^{m-3}\leq(q-1)(q-b+1)q^{m-2}$$ which means that $k\geq q-b+2$. These $(q-b+2)$ hyperplanes are minimal weight codewords of $R_q(q-1+b,m)$. So, they meet $S$ in $(q-b)$ affine subspaces of codimension 2 parallel to $Q$, that is to say, they meet $S\cap G_1$ in $(q-b)q^{m-3}$ points. This is absurd since $S\cap G_1$ is of type $B$ and so there are at most $(q-b+1)$ affine hyperplanes through $Q$ which meet $S\cap G_1$ in $(q-b)q^{m-3}$ points (see Figure \ref{fig8}). So there exists an affine hyperplane through $Q$ which does not meet $S$. 
\begin{figure}[!h]
\caption{}
\label{fig8}
\begin{center}\begin{tikzpicture}[scale=0.2]
\draw (0,0)--(5,2)--(5,11)--(0,9)--cycle;
\draw (1,2/5)--(1,9+2/5);
\draw (2,4/5)--(2,9+4/5);
\draw[dotted] (3,6/5)--(3,9+6/5);
\draw[dotted] (4,8/5)--(4,9+8/5);
\draw (0,9) node[above left]{$H_1$};
\draw[dashed](0,4)--++(5,2);
\draw (4,9+8/5) node[above]{$Q$};
\draw (0,4)--(12,4)--++(5,2)--++(-12,0);
\draw (17,6) node[above right]{$G_{1}$};
\draw (2,4+4/5)--++(12,0);
\draw (1,4+2/5)--(15,4+6/5);
\draw[dotted] (3,4+6/5)--(13,4+2/5);
\draw[dotted] (4,4+8/5)--(12,4);
\draw[dotted] (6,4)--(11,6);
\draw (6,0)--++(5,2)--++(0,9)--++(-5,-2)--cycle;
\draw (11,11) node[above]{$H$};
\end{tikzpicture}\end{center}
\end{figure}

By applying the same argument to all affine subspaces of codimension 2 included in $H_1$ which does not meet $S$, since $S\cap G_i$ is of type $B$ for all $i$, we get that $S$ is of type $B$.
\end{itemize}
\end{preuve}

\subsection{The support is included in an affine subspace of codimension $t$.}

The two following lemmas are proved in \cite{erickson1974counting}.

\begin{lemme}\label{3.5}Let $m\geq2$, $q\geq3$, $1\leq t\leq m-1$, $1\leq s\leq q-2$. Assume that $f\in R_q(t(q-1)+s,m)$ is such that $\forall x=(x_1,\ldots,x_m)\in\F_q^m$, $$f(x)=(1-x_1^{q-1})\widetilde{f}(x_2,\ldots,x_m)$$ and that $g\in R_q(t(q-1)+s-k)$, $1\leq k\leq q-1$, is such that $(1-x_1^{q-1})$ does not divide $g$. Then, if $h=f+g$, either $|h|\geq (q-s+k)q^{m-t-1}$ or $k=1$.\end{lemme}

\begin{lemme}\label{3.6}Let $m\geq2$, $q\geq3$, $1\leq t\leq m-1$, $1\leq s\leq q-2$ and \\$f\in R_q(t(q-1)+s,m)$. For $a\in\F_q$, the function $f_a$ of $B_{m-1}^q$ defined for all $(x_2,\ldots,x_m)\in\F_q^m$ by $f_a(x_2,\ldots,x_m)=f(a,x_2,\ldots,x_m)$. Assume that for $a$, $b\in\F_q$ $f_a$ is different from the zero function and $(1-x_2^{q-1})$ divides $f_a$ and that $$0<|f_b|<(q-s+1)q^{m-t-2}.$$ Then there exists $T$ an affine transformation, fixing $x_i$ for $i\neq2$ such that $(1-x_2^{q-1})$ divides $(f\circ T)_a$ and $(f\circ T)_b$.\end{lemme}

\begin{lemme}\label{hyp}Let $m\geq3$, $q\geq4$, $1\leq t\leq m-2$ and $2\leq s \leq q-2$. If $f\in R_q(t(q-1)+s,m)$ is such that $|f|=(q-s+1)(q-1)q^{m-t-2}$, then the support of $f$ is included in an affine hyperplane of $\F_q^m$.\end{lemme}

\begin{preuve} We denote by $S$ the support of $f$. Assume that $S$ is not included in an affine hyperplane. Then, by Lemma \ref{inter}, there exists an affine hyperplane $H$ such that either $H$ does not meet $S$ or $H$ meets $S$ in $(q-s)q^{m-t-2}$. Now, by Lemma \ref{cle}, since $S$ is not included in an affine hyperplane, either $S$ meets all affine hyperplanes parallel to $H$ or $S$ meets $(q-1)$ affine hyperplanes parallel to $H$ in $(q-s+1)q^{m-t-2}$ or $S$ meets $(q-s+1)$ affine hyperplanes parallel to $H$ in $(q-1)q^{m-t-2}$ points. By applying an affine transformation, we can assume that $x_1=\lambda$, $\lambda \in \F_q$ is an equation of $H$. We define $f_{\lambda}\in B_{m-1}^q$ by $$\forall(x_2,\ldots,x_m)\in\F_q^{m-1}, \qquad f_{\lambda}(x_2,\ldots,x_m)=f(\lambda,x_2,\ldots,x_m).$$ We set an order $\lambda_1,\ldots,\lambda_q$ on the elements of $\F_q$ such that $$|f_{\lambda_1}|\leq\ldots\leq|f_{\lambda_q}|.$$
Then either $|f_{\lambda_1}|=0$ or $|f_{\lambda_1}|=(q-s)q^{m-t-2}$, that is to say either $f_{\lambda_1}$ is null or $f_{\lambda_1}$ is the minimal weight codeword of $R_q(t(q-1)+s,m-1)$ and its support is included in an affine subspace of codimension $t+1$. Since $t\geq1$, in both cases, the support of $f_{\lambda_1}$ is included in an affine hyperplane of $\F_q^m$ different from the hyperplane parallel to $H$ of equation $x_1=\lambda_1$. By applying an affine transformation that fixes $x_1$, we can assume that $(1-x_2^{q-1})$ divides $f_{\lambda_1}$. Since $S$ is not included in an affine hyperplane, there exists $2\leq k\leq q$ such that $1-x_2^{q-1}$ does not divide $f_{\lambda_k}$. We denote by $k_0$ the smallest such $k$. 

Assume that $S$ meets all affine hyperplanes parallel to $H$ and  that $$|f_{\lambda_{k_0}}|\geq(q-s+k_0-1)q^{m-t-2}.$$ Then 
\begin{align*}|f|&=\sum_{k=1}^q|f_{\lambda_k}|\\&\geq(q-s)q^{m-t-2}(k_0-1)+(q-k_0+1)(q-s+k_0-1)q^{m-t-2}\\&=(q-s)q^{m-t-1}+(k_0-1)(q-k_0+1)q^{m-t-2}\\&>(q-s)q^{m-t-1}+(s-1)q^{m-t-2}\end{align*}
which gives a contradiction. In the cases where $S$ meets $(q-s')$, $s'=1$ or $s'=s-1$, for $1\leq i\leq s'$, $|f_{\lambda_i}|=0$ and the support of $f_{\lambda_{s'+1}}$ is $S\cap H_{\lambda_{s'+1}}$, where $H_{\lambda_{s'+1}}$ is the hyperplane of equation $x_1=\lambda_{s'+1}$. Since $S\cap H_{\lambda_{s'+1}}$ is the support of a minimum weight codeword of $R_q((t+1)(q-1)+s',m)$, it is included in affine subspace of codimension $t+1$. So in those cases, we can assume that $k_0\geq s'+2$.  Finally, $|f_{\lambda_{k_0}}|<(q-s+k_0-1)q^{m-t-2}$.

We write \begin{align*}f(x_1,x_2,x_3,\ldots,x_m)&=\sum_{i=0}^{q-1}x_2^ig_i(x_1,x_3,\ldots,x_m)
\\&=h(x_1,x_2,x_3,\ldots,x_m)+(1-x_2^{q-1})g(x_1,x_3,\ldots,x_m).\end{align*}
Since for all $1\leq i\leq k_0-1$, $1-x_2^{q-1}$ divides $f_{\lambda_i}$, for all $(x_2,\ldots,x_m)\in\F_q^{m-1}$, for all $1\leq i\leq k_0-1$, $h(\lambda_i,x_2,\ldots,x_m)=0$. So, by Lemma \ref{DGMW1}, \begin{align*}f(x_1,x_2,x_3,\ldots,x_m)&=(x_1-\lambda_1)\ldots(x_1-\lambda_{k_0-1})\widetilde{h}(x_1,x_2,x_3,\ldots,x_m)\\&\hspace{2cm}+(1-x_2^{q-1})g(x_1,x_3,\ldots,x_m)\end{align*} with $\deg(\widetilde{h})\leq r-k_0+1$. Then by applying Lemma \ref{3.5} to $f_{\lambda_{k_0}}$, since $$|f_{\lambda_{k_0}}|<(q-s+k_0-1)q^{m-t-2},$$ $k_0=2$. This gives a contradiction in the cases where $S$ does not meet all hyperplanes parallel to $H$. In the case where $S$ meets all hyperplanes parallel to $H$, by applying Lemma \ref{3.6}, there exists $T$ an affine transformation which fixes $x_1$ such that $(1-x_2^{q-1})$ divides $(f\circ T)_{\lambda_1}$ and $(f\circ T)_{\lambda_2}$, we set $k_0'$ the smallest $k$ such that $(1-x_2^{q-1})$ does not divide $(f\circ T)_{\lambda_k}$. Then $k_0'\geq3$ and by applying the previous argument to $f \circ T$, we get a contradiction. \end{preuve}

\begin{proposition}\label{inclu}Let $m\geq3$, $q\geq4$, $1\leq t\leq m-2$ and $2\leq s \leq q-2$. If $f\in R_q(t(q-1)+s,m)$ is such that $|f|=(q-1)(q-s+1)q^{m-t-2}$, then the support of $f$ is included in an affine subspace of codimension $t$.\end{proposition}

\begin{preuve}We denote by $S$ the support of $f$. By Lemma \ref{hyp}, $S$ is included in $H$ an affine hyperplane. By applying an affine transformation, we can assume that $x_1=0$ is an equation of $H$. Let $g\in B_{m-1}^q$ defined by $$\textrm{$\forall x=(x_{2},\ldots,x_m)\in\F_{q}^{m-1}$, $g(x)=f(0,x_{2},\ldots,x_m)$}$$ and denote by $P\in\F_{q}[X_{2},\ldots,X_m]$ its reduced form. Since $$\forall x=(x_1,\ldots,x_m)\in\F_{q}^m, \ f(x)=(1-x_1^{q-1})P(x_{2},\ldots,x_m),$$ the reduced form of $f\in R_q(t(q-1)+s,m)$ is $$(1-X_1^{q-1})P(X_{2},\ldots,X_m).$$ Then $g\in R_q((t-1)(q-1)+s,m-1)$ and $$|g|=|f|=(q-s+1)(q-1)q^{m-t-2}=(q-1)(q-s+1)q^{m-1-(t-1)-2}.$$ Then, by Lemma \ref{hyp}, if $t\geq2$, the support of $g$ is included in an affine hyperplane of $\F_q^{m-1}$. By iterating this argument, we get that $S$ is included in an affine subspace of codimension $t$. \end{preuve}

\subsection{Proof of Theorem \ref{wpoids2}}

Let $0\leq t\leq m-2$, $2\leq s\leq q-2$ and $f\in R_q(t(q-1)+s,m)$ such that $$|f|=(q-s+1)(q-1)q^{m-t-2};$$ we denote by $S$ the support of $f$. Assume that $t\geq1$. By Proposition \ref{inclu}, $S$ is included in an affine subspace $G$ of codimension $t$. By applying an affine transformation, we can assume that $$G=\{x=(x_1,\ldots,x_m)\in\F_{q}^m:x_{i}=0 \textrm{ for }1\leq i\leq t\}.$$ Let $g\in B_{m-t}^q$ defined for all $x=(x_{t+1},\ldots,x_m)\in\F_{q}^{m-t}$ by $$g(x)=f(0,\ldots,0,x_{t+1},\ldots,x_m)$$ and denote by $P\in\F_{q}[X_{t+1},\ldots,X_m]$ its reduced form. Since $$\forall x=(x_1,\ldots,x_m)\in\F_{q}^m, \ f(x)=(1-x_1^{q-1})\ldots(1-x_t^{q-1})P(x_{t+1},\ldots,x_m),$$ the reduced form of $f\in R_q(t(q-1)+s,m)$ is $$(1-X_1^{q-1})\ldots(1-X_t^{q-1})P(X_{t+1},\ldots,X_m).$$ Then $g\in R_q(s,m-t)$ and $|g|=|f|=(q-s+1)(q-1)q^{m-t-2}$. Thus, using the case where $t=0$, we finish the proof of Theorem \ref{wpoids2}.

\section{Case where $s=0$}

\subsection{The support is included in an affine subspace of dimension $m-t+1$}

\begin{proposition}\label{m-10}Let $q\geq 3$, $m\geq2$ and $f\in R_q((m-1)(q-1),m)$ such that $|f|=2(q-1)$. Then, the support of $f$ is included in an affine plane.\end{proposition}

In order to prove this proposition, we need the following lemma.

\begin{lemme}\label{lam-1}Let $m\geq3$, $q\geq4$ and $f\in R_q((m-1)(q-1),m)$ such that $|f|=2(q-1)$. If $H$ is an affine hyperplane of $\F_q^m$ such that $S\cap H\neq S$, $\#(S\cap H)=N$, $3\leq N\leq q-1$ and $S\cap H$ is not included in a line then there exists $H_1$ an affine hyperplane of $\F_q^m$ such that $S\cap H_1\neq S$, $\#(S\cap H_1)\geq N+1$ and $S\cap H_1$ is not included in a line\end{lemme}

\begin{preuve}Since $S\cap H\neq S$, by Lemma \ref{cle}, either $S$ meets $(q-1)$ hyperplanes parallel to $H$ or $S$ meets 2 hyperplanes parallel to $H$ or $S$ meets all affine hyperplanes parallel to $H$. If $S$ does not meet all affine hyperplanes parallel to $H$ then $S\cap H$ is the support of a minimal weight codeword of \\$R_q((m-1)(q-1)+s',m)$, $s'=1$ or $s'=q-2$. In both cases, $S\cap H$ is included in a line which is absurd. So, $S$ meets all affine hyperplanes parallel to $H$. 

By applying an affine transformation, we can assume that $x_1=0$ is an equation of $H$. Let $I:=\{a\in\F_q:\#(\{x_1=a\}\cap S)=1\}$ and $k:=\#I$. Since $\#S=2(q-1)$ and $\#(S\cap H)=N$, $k\geq N$. We define $$\forall x=(x_1,\ldots,x_m)\in\F_q^m,\quad g(x)=f(x)\prod_{a\not\in I}(x_1-a).$$
Then, $\deg(g)\leq (m-1)(q-1)+q-k$ and $|g|=k$. So, $g$ is a minimal weight codewords of $R_q((m-1)(q-1)+q-k,m)$ and its support is included in a line $L$ which is not included in $H$. We denote by $\overrightarrow{u}$  a directing vector of $L$. Let $b$ be the intersection point of $H$ and $L$ and $\omega_1$, $\omega_2$, $\omega_3$ 3 points of $S\cap H$ which are not included in a line. Then there exists $\overrightarrow{v}$ and $\overrightarrow{w}\in\{\overrightarrow{b\omega_1},\overrightarrow{b\omega_2},\overrightarrow{b\omega_3}\}$ which are linearly independent. Since $L$ is not included in $H$, $\{\overrightarrow{u},\overrightarrow{v},\overrightarrow{w}\}$ are linearly independent. We choose $H_1$ an affine hyperplane such that $b\in H_1$, $b+\overrightarrow{v}\in H_1$, $L\subset H_1$ but $b+\overrightarrow{w}\not \in H_1$.  \end{preuve}

Now we can prove the proposition

\begin{preuve}If $m=2$, we have the result. Assume $m\geq3$. Let $S$ be the support of $f$. Since $\#S=2(q-1)>q$, $S$ is not included in a line. Let $\omega_1$, $\omega_2$, $\omega_3$ 3 points of $S$ not included in a line. Let $H$ be an hyperplane such that $\omega_1$, $\omega_2$, $\omega_3\in H$. Assume that $S\cap H\neq S$. Then there exists an affine hyperplane $H_1$ such that $\#(S\cap H_1)\geq q$, $S\cap H_1$ is not included in a line and $S\cap H_1\neq S$. Indeed, if $q=3$, we take $H_1=H$ and for $q\geq4$, we proceed by induction using the previous Lemma. Then by Lemma \ref{cle} either $S$ meets 2 hyperplanes parallel to $H_1$ in 2 points or $S$ meets 2 hyperplanes parallel to $H_1$ in $q-1$ points or $S$ meets all affine hyperplanes parallel to $H_1$. Since $\#(S\cap H_1)\geq q$, $S$ meets all hyperplanes parallel to $H_1$. Then, we must have $$q+q-1\leq2(q-1)$$ which is absurd.  \end{preuve}

The two following lemmas are proved in \cite{erickson1974counting}.

\begin{lemme}\label{2.14}Let $m\geq2$, $q\geq3$, $1\leq t\leq m$ and $f\in R_q(t(q-1),m)$ such that $|f|=q^{m-t}$ and $g\in R_q((t(q-1)-k,m)$, $1\leq k\leq q-1$, such that $g\neq0$. If $h=f+g$ then either $|h|=kq^{m-t}$ or $|h|\geq(k+1)q^{m-t}$.\end{lemme}

\begin{lemme}\label{2.15}Let $m\geq2$, $q\geq3$, $1\leq t\leq m-1$ and $f\in R_q(t(q-1),m)$. For $a\in\F_q$, we define the function $f_a$ of $B_{m-1}^q$ by for all $(x_2,\ldots,x_m)\in\F_q^m$, $f_a(x_2,\ldots,x_m)=f(a,x_2,\ldots,x_m)$. If for some $a$, $b\in\F_q$, $|f_a|=|f_b|=q^{m-t-1}$, then there exists $T$ an affine transformation fixing $x_1$ such that $$(f\circ T)_a=(f\circ T)_b.$$\end{lemme}

\begin{proposition}\label{inclu0}Let $q\geq3$, $m\geq2$, $1\leq t\leq m-1$. If $f\in R_q(t(q-1),m)$ is such that $|f|=2(q-1)q^{m-t-1}$ then the support of $f$ is included in an affine subspace of dimension $m-t+1$. \end{proposition}

\begin{preuve}For $t=1$, this is obvious. For the other cases we proceed by recursion on $t$. Proposition \ref{m-10} gives the case where $t=m-1$.

If $m\leq 3$ we have considered all cases. Assume $m\geq4$. Let $2\leq t\leq m-2$. Assume that for $f\in R_q((t+1)(q-1),m)$ such that $|f|=2(q-1)q^{m-t-2}$ the support of $f$ is included in an affine subspace of dimension $m-t$. Let $f\in R_q(t(q-1),m)$ such that $|f|=2(q-1)q^{m-t-1}$. We denote by $S$ the support of $f$. 

Assume that $S$ is not included in an affine subspace of dimension $m-t+1$. Then there exists $H$ an affine hyperplane of $\F_q^m$ such that $S\cap H\neq S$ and $S\cap H$ is not included in an affine space of dimension $m-t$. By Lemma \ref{cle}, either $S$ meets all affine hyperplanes parallel to $H$ or $S$ meets $(q-1)$ affine hyperplanes parallel to $H$ in $2q^{m-t-1}$ or $S$ meets 2 affine hyperplanes parallel to $H$ in $(q-1)q^{m-t-1}$ points. 

If $S$ does not meet all hyperplanes parallel to $H$ then $S\cap H$ is the support of a minimal weight codeword of $R_q(t(q-1)+s',m)$, $s'=1$ or $s'=q-2$. So $S\cap H$ is included in an affine subspace of dimension $m-t$ which gives a contradiction. 

So, $S$ meets all affine hyperplanes parallel to $H$ in at least $q^{m-t-1}$ points. If for all $H'$ parallel to $H$, $\#(S\cap H')>q^{m-t-1}$ then for all $H'$ parallel to $H$, $\#(S\cap H')\geq2(q-1)q^{m-t-2}$. So, for reason of cardinality, $S\cap H$ is the support of a second weight codeword of $R_q((t+1)(q-1),m)$ and by recursion hypothesis $S\cap H$ is included in an affine subspace of dimension $m-t$ which gives a contradiction. So, there exists $H_1$ an affine hyperplane parallel to $H$ such that $\#(S\cap H_1)=q^{m-t-1}$.

By applying an affine transformation, we can assume that $x_1=\lambda$, $\lambda\in \F_q$ is an equation of $H$. For $\lambda\in\F_q$, we define $f_{\lambda}\in B_{m-1}^q$ by $$\forall (x_2,\ldots,x_m)\in \F_q^{m-1},\qquad f_{\lambda}(x_2,\ldots,x_m)=f(\lambda,x_2,\ldots,x_m).$$ 
We set an order $\lambda_1,\ldots, \lambda_q$ on the elements of $\F_q$ such that 
$$|f_{\lambda_1}|\leq\ldots\leq|f_{\lambda_q}|.$$ 
Since $\#(S\cap H_1)=q^{m-t-1}$ and $S$ meets all hyperplanes parallel to $H$, $$|f_{\lambda_1}|=q^{m-t-1}$$ and $f_{\lambda_1}$ is a minimum weight codeword of $R_q(t(q-1),m-1)$. Let $k_0$ be the smallest integer such that $|f_{\lambda_{k_0}}|>q^{m-t-1}$. Since $|f|>q^{m-t}$, $k_0\leq q$. Then by Lemma \ref{2.15} and applying an affine transformation that fixes $x_1$, we can assume that for all $2\leq i\leq k_0-1$, $f_{\lambda_i}=f_{\lambda_1}$. If we write for all $x=(x_1,\ldots,x_m)\in\F_q^m$, $$f(x)=f_{\lambda_1}(x_2,\ldots,x_m)+(x_1-\lambda_1)\widehat{f}(x_1,\ldots,x_m).$$
Then for all $2\leq i\leq k_0-1$, for all $\overline{x}=(x_2,\ldots,x_m)\in\F_q^{m-1}$, $$f_{\lambda_i}(\overline{x})=f_{\lambda_1}(\overline{x})+(\lambda_i-\lambda_1)\widehat{f}_{\lambda_i}(\overline{x}).$$
Since for all $2\leq i\leq k_0-1$, $f_{\lambda_i}=f_{\lambda_1}$, by Lemma \ref{DGMW1}, we can write for all 
$ x=(x_1,\ldots,x_m)\in\F_q^m$,$$f(x)=f_{\lambda_1}(x_2,\ldots,x_m)+(x_1-\lambda_1)\ldots(x_1-\lambda_{k_0-1})\overline{f}(x_1,\ldots,x_m)$$ with $\deg(\overline{f})\leq t(q-1)-k_0+1$. Now, we have $f_{\lambda_{k_0}}=f_{\lambda_1}+\lambda'\overline{f}_{\lambda_{k_0}}$, $\lambda'\in\F_q^*$. Then, by Lemma \ref{2.14}, either $|f_{\lambda_{k_0}}|\geq k_0q^{m-t-1}$ or $|f_{\lambda_{k_0}}|=(k_0-1)q^{m-t-1}$. Assume that $|f_{\lambda_{k_0}}|\geq k_0q^{m-t-1}$. Then \begin{align*}|f|&=\sum_{i=1}^q|f_{\lambda_i}|\\&\geq (k_0-1)q^{m-t-1}+(q+1-k_0)k_0q^{m-t-1}\\&=q^{m-t}+(k_0-1)(q-k_0+1)q^{m-t-1}\\&>2(q-1)q^{m-t-1}.\end{align*}
So, $|f_{\lambda_{k_0}}|=(k_0-1)q^{m-t-1}$. Since $|f_{\lambda_{k_0}}|>q^{m-t-1}$, $k_0\geq3$. Now, we have $$|f|\geq (k_0-1)q^{m-t-1}+(q+1-k_0)(k_0-1)q^{m-t-1}=(k_0-1)(q-k_0+2)q^{m-t-1}.$$ So either $k_0=q$ or $k_0=3$.
\begin{itemize}\item Assume $k_0=q$.\\
Since $f_{\lambda_1}=\ldots=f_{\lambda_{q-1}}$ are minimum weight codeword of \\$R_q(t(q-1),m-1)$, there exists $A$ an affine subspace of dimension $m-t$ of $\F_q^m$ such that for all $1\leq i\leq q-1$, $S\cap H_i\subset A$, where $H_i$ is the hyperplane parallel to $H$ of equation $x_1=\lambda_i$. Since $S$ is not included in an affine subspace of dimension $m-t+1$ and $t\geq2$, there exists an affine hyperplane $G$ containing $A$ such that $S\cap G\neq S$ and there exists $x\in S\cap G$, $x\not\in A$. Then $\#(S\cap G)\geq (q-1)q^{m-t-1}+1$, $S\cap G\neq S$ and $S\cap G$ is not included in an affine subspace of dimension $m-t$. Applying to $G$ the same argument than to $H$, we get a contradiction.
\item So, $k_0=3$.\\
Then $f_{\lambda_1}=f_{\lambda_2}$ are minimum weight codeword of $R_q(t(q-1),m-1)$ and for reason of cardinality, for all $3\leq i\leq q$, $|f_{\lambda_i}|=2q^{m-t-1}$. So, there exists $A$ an affine subspace of dimension $m-t$ of $\F_q^m$ such that for all $1\leq i\leq 2$, $S\cap H_i\subset A$, where $H_i$ is the hyperplane parallel to $H$ of equation $x_1=\lambda_i$. Since $S$ is not included in an affine subspace of dimension $m-t+1$ and $t\geq2$, there exists an affine hyperplane $G$ containing $A$ such that $S\cap G\neq S$ and there exists $x\in S\cap G$, $x\not\in A$. Then $\#(S\cap G)\geq 2q^{m-t-1}+1$, $S\cap G\neq S$ and $S\cap G$ is not included in an affine subspace of dimension $m-t$. Applying to $G$ the same argument than to $H$, we get a contradiction.
\end{itemize} 
Finally, $S$ is included in an affine subspace of dimension $m-t+1$.\end{preuve}

\subsection{Proof of Theorem \ref{0}}

Let $1\leq t\leq m-1$ and $f\in R_q(t(q-1),m)$ such that $$|f|=2(q-1)q^{m-t-1};$$ we denote by $S$ the support of $f$. Assume that $t\geq2$. By proposition \ref{inclu0}, $S$ is included in an affine subspace $G$ of codimension $t-1$. By applying an affine transformation, we can assume that $$G=\{x=(x_1,\ldots,x_m)\in\F_{q}^m:x_{i}=0 \textrm{ for }1\leq i\leq t-1\}.$$ Let $g\in B_{m-t+1}^q$ defined for all $x=(x_{t},\ldots,x_m)\in\F_{q}^{m-t+1}$ by $$g(x)=f(0,\ldots,0,x_{t},\ldots,x_m)$$ and denote by $P\in\F_{q}[X_{t},\ldots,X_m]$ its reduced form. Since $$\forall x=(x_1,\ldots,x_m)\in\F_{q}^m, \ f(x)=(1-x_1^{q-1})\ldots(1-x_{t-1}^{q-1})P(x_{t},\ldots,x_m),$$ the reduced form of $f\in R_q(t(q-1)+s,m)$ is $$(1-X_1^{q-1})\ldots(1-X_{t-1}^{q-1})P(X_{t},\ldots,X_m).$$ Then $g\in R_q(q-1,m-t+1)$ and $|g|=|f|=2(q-1)q^{m-t-1}$. Thus, using the case where $t=1$, we finish the proof of Theorem \ref{0}.

\section{Case where $0\leq t\leq m-2$ and $s=1$}

\subsection{Case where $q\geq4$}

\begin{lemme}\label{cle1}Let $m\geq2$, $q\geq4$, $0\leq t\leq m-2$ and $f\in R_q(t(q-1)+1,m)$ such that $|f|=q^{m-t}$. We denote by $S$ the support of $f$. Then, if $H$ is an affine hyperplane of $\F_q^m$ such that $S\cap H\neq \emptyset$ and $S\cap H\neq S$, $S$ meets all affine hyperplanes parallel to $H$.\end{lemme}

\begin{preuve}By applying an affine transformation, we can assume that $x_1=0$ is an equation of $H$. Let $H_a$ be the $q$ affine hyperplanes parallel to $H$ of equation $x_1=a$, $a\in\F_q$. We denote by $I:=\{a\in\F_q:S\cap H_a=\emptyset\}$. Let $k:=\#I$ and assume that $k\geq1$. Since $S\cap H\neq\emptyset$ and $S\cap H\neq S$, $k\leq q-2$. For all $c\not\in I$ we define $$\forall x=(x_1,\ldots,x_m)\in\F_q^m, \quad g_c(x)=f(x)\prod_{a\in\F_q\setminus I, a\neq c}(x_1-a).$$
Then $|f|=\displaystyle\sum_{c\not\in I}|g_c|$.
\begin{itemize}\item Assume $k\geq2$.
\\Then for all $c\not\in I$, $\deg(g_c)\leq t(q-1)+q-k$ and $2\leq q-k\leq q-2$. So, $|g_c|\geq kq^{m-t-1}$. Let $N=\#\{c\not\in I:|g_c|=kq^{m-t-1}\}$. If $|g_c|>kq^{m-t-1}$, $|g_c|\geq (k+1)(q-1)q^{m-t-2}$. Hence
$$q^{m-t}\geq Nkq^{m-t-1}+(q-k-N)(k+1)(q-1)q^{m-t-2}.$$
Since $k\geq2$, we get that $N\geq q-k$. Since $(q-k)kq^{m-t-1}\neq q^{m-t}$, we get a contradiction.
\item Assume $k=1$.
\\Then, for all $c\not \in I$, $\deg(g_c)\leq t(q-1)+1+q-2=(t+1)(q-1)$. So $|g_c|\geq q^{m-t-1}$. Let $N=\#\{c\not\in I:|g_c|=q^{m-t-1}\}$. If $|g_c|>q^{m-t-1}$, $|g_c|\geq 2(q-1)q^{m-t-2}$. Since for $q\geq4$, $2(q-1)^2q^{m-t-2}>q^{m-t}$, $N\geq1$. Furthermore, since $(q-1)q^{m-t-1}<q^{m-t}$, $N\leq q-2$. For $\lambda\in\F_q$, we define $f_{\lambda}\in B_{m-1}^q$ by $$\forall (x_2,\ldots,x_m)\in \F_q^{m-1},\qquad f_{\lambda}(x_2,\ldots,x_m)=f(\lambda,x_2,\ldots,x_m).$$ We set $\lambda_1,\ldots,\lambda_q$ an order on the elements of $\F_q$ such that for all $i\leq N$, $|f_{\lambda_i}|=q^{m-t-1}$, $|f_{\lambda_{N+1}}|=0$ and $q^{m-t-1}<|f_{\lambda_{N+2}}|\leq \ldots \leq|f_{\lambda_q}|$.

Since $f_{\lambda_{N+1}}=0$, we can write for all $(x_1,\ldots,x_m)\in\F_q^m$, $$f(x_1,\ldots,x_m)=(x_1-\lambda_{N+1})h(x_1,\ldots,x_m)$$ with $\deg(h)\leq t(q-1)$. Then, for all $1\leq i\leq q$, $i\neq N+1$ and $(x_2,\ldots,x_m)\in\F_q^{m-1}$, $$f_{\lambda_i}(x_2,\ldots,x_m)=(\lambda_i-\lambda_{N+1})h_{\lambda_i}(x_2,\ldots,x_m).$$ So $\deg(f_{\lambda_i})\leq t(q-1)$ and $h_{\lambda_i}=\frac{f_{\lambda_i}}{\lambda_i-\lambda_{N+1}}$. 

Since $h\in R_q(t(q-1),m)$, by Lemma \ref{2.15}, there exists an affine transformation such that for all $i\leq N$, $h_{\lambda_i}=h_{\lambda_1}$. Then, for all $(x_1,\ldots,x_m)\in\F_q^m$, $$h(x_1,\ldots,x_m)=h_{\lambda_1}(x_2,\ldots,x_m)+(x_1-\lambda_1)\ldots(x_1-\lambda_N)\widetilde{h}(x_1,\ldots,x_m)$$ with $\deg(\widetilde{h})\leq t(q-1)-N$. Hence, for all $(x_1,\ldots,x_m)\in\F_q^m$, $$f(x_1,\ldots,x_m)=\frac{x_1-\lambda_{N+1}}{\lambda_1-\lambda_{N+1}}f_{\lambda_1}(x_2\ldots,x_m)+(x_1-\lambda_1)\ldots(x_1-\lambda_{N+1})\widetilde{h}(x_1,\ldots,x_m).$$
Then, for all $(x_2,\ldots,x_m)\in\F_q^{m-1}$, $$f_{\lambda_{N+2}}(x_2,\ldots,x_m)=\lambda f_{\lambda_1}(x_2\ldots,x_m)+\lambda'\widetilde{h}_{\lambda_{n+2}}(x_2,\ldots,x_m)$$ with $\lambda$, $\lambda'\in\F_q^*$. 

Since $f_{\lambda_1}\in R_q(t(q-1),m-1)$ and $\widetilde{h}_{\lambda_{n+2}}\in R_q(t(q-1)-N,m-1)$, by Lemma \ref{2.14}, either $|f_{\lambda_{N+2}}|=Nq^{m-t-1}$ or $|f_{\lambda_{N+2}}|\geq (N+1)q^{m-t-1}$.

If $N=1$, since $|f_{\lambda_{N+2}}|>q^{m-t-1}$, we get $$q^{m-t-1}+(q-2)2q^{m-t-1}\leq q^{m-t}$$ which means that $q\leq 3$. So $N\geq 2$. Then, $$Nq^{m-t-1}+(q-1-N)Nq^{m-t-1}\leq q^{m-t}.$$ Since $N(q-N)\geq 2(q-2)$, we get that $q\leq 4$. So, the only possibility is $q=4$ and $N=q-2=2$.

If $t=0$, $H_{\lambda_4}$ contains $2.4^{m-1}$ points which is absurd. Assume $t\geq1$. Since $h_{\lambda_1}=h_{\lambda_2}$ and for $i\in\{1,2\}$, $f_{\lambda_i}=(\lambda_i-\lambda_3)h_{\lambda_i}$, $S\cap H_{\lambda_1}$ and $S\cap H_{\lambda_2}$ are both included in $A$ an affine subspace of dimension $m-t$. If $t=1$, by applying an affine transformation which fixes $x_1$, we can assume that \\$x_2=0$ is an equation of $A$. If $S$ is included in $A$, then $$\#(S\cap H_{\lambda_4}\cap A)=2.4^{m-2}$$ which is absurd since $H_{\lambda_4}\cap A$ is an affine subspace of codimension 2. So there exists an affine hyperplane $H'$ containing $A$ but not $S$. By applying an affine transformation which fixes $x_1$, we can assume that $x_2=0$ is an equation of $H'$. Now, consider $g$ defined for all $(x_1,\ldots,x_m)\in\F_q^m$ by $g(x_1,\ldots,x_m)=x_2f(x_1,\ldots,x_m)$. Then $|g|\leq 2.4^{m-t-1}$. Furthermore, since $S$ is not included in $H'$ and $\deg(g)\leq 3t+2$, $|g|\geq 2.4^{m-t-1}$. So $g$ is a minimum weight codeword of $R_4(3t+2,m)$ and its support is the union of 2 parallel affine subspace of codimension $t+1$ included in an affine subspace of codimension $t$. Then, since $H'\cap H_{\lambda_4}=\emptyset$, there exists $H'_1$ an hyperplane parallel to $H'$ such that $S\cap H'_1=\emptyset$. Now, consider $G$ the hyperplane through $H_{\lambda_4}\cap H'_1$ and $H'\cap H_{\lambda_3}$ and $G'$ the hyperplane through $H'\cap H_{\lambda_4}$ parallel to $G$ (see Figure \ref{fig9}).
\begin{figure}[!h]
\caption{}
\label{fig9}
\begin{center}\begin{tikzpicture}[scale=0.2]
\draw (2,1/2)--(10,4+1/2)--++(0,9)--(2,9+1/2)--cycle;
\draw (2,1/2) node[below]{$H'$};
\draw (3,1)--++(0,9);
\draw (5,2)--++(0,9);
\draw[dotted] (7,3)--++(0,9);
\draw[dotted] (9,4)--++(0,9);
\draw (16,1/2)--(8,4+1/2)--++(0,9)--(16,9+1/2)--cycle;
\draw (16,1/2) node[below right]{$H_{\lambda_4}$};
\draw (15,1)--++(0,9);
\draw (13,2)--++(0,9);
\draw[dotted] (11,3)--++(0,9);
\draw (0,3)--++(19,0)--++(0,9)--++(-19,0)--cycle;
\draw (19,12) node[right]{$G$};
\draw (1,4)--++(17,0)--++(0,9)--++(-17,0)--cycle;
\draw (18,13) node[above]{$G'$};
\end{tikzpicture}\end{center}
\end{figure}

 Then $G$ and $G'$ does not meet $S$ but $S$ is not included in an hyperplane parallel to $G$ which is absurd by the previous case.
\end{itemize}\end{preuve}
\pagebreak
\begin{lemme}\label{m-2}For $m\geq3$, if $f\in R_4(3(m-2)+1,m)$ is such that $|f|=16$, the support of $f$ is an affine plane.\end{lemme}

\begin{preuve}We denote by $S$ the support of $f$.

First, we prove the case where $m=3$. To prove this case, by Lemma \ref{cle1}, we only have to prove that there exists an affine hyperplane which does not meet $S$. 

Assume that $S$ meets all affine hyperplanes. Let $H$ be an affine hyperplane. Then for all $H'$ affine hyperplane parallel to $H$, $\#(S\cap H')\geq3$. Assume that for all $H'$ hyperplane parallel to $H$, $\#(S\cap H')\geq4$. For reason of cardinality , for all $H'$ parallel to $H$ $\#(S\cap H')=4$. Since $q=4$, there exists a line in $H$ which does not meet $S$. Since $3.4+4=16$, $S$ meets 4 hyperplanes through this line in 3 points and the last one in 4 points. So, there exists $H_1$ an affine hyperplane such that $\#(S\cap H_1)=3$. We denote by $H_2$, $H_3$, $H_4$ the hyperplanes parallel to $H_1$.
 Then, $S\cap H_1$ is the support of a minimal weight codeword of $R_4(3(m-1)+1,m)$ so $S\cap H_1$ is included in $L$ a line. Consider $L'$ a line in $H_1$ parallel to $L$. Then there is 4 hyperplanes through $L'$ which meets $S$ in 3 points and one $H_1'$ which meets $S$ in 4 points. Let $H'$ be an affine hyperplane through $L'$ which meets $S$ in 3 points; $S\cap H'$ is minimum weight codeword of $R_4(3(m-1)+1,m)$ which does not meet $H_1$. So either $S\cap H'$ is included in an affine hyperplane parallel to $H_1$ or $S\cap H'$ meets all affine hyperplane parallel to $H_1$ but $H_1$ in 1 point. Then we consider 4 cases :
\begin{enumerate}\item $H_1$ is the only hyperplane through $L'$ such that $\#(S\cap H_1)=3$ and $S\cap H_1$ is included in one of the affine hyperplane parallel to $H_1$.
\\Since $S\cap H_1\cap H_1'=\emptyset$ there exists an affine hyperplane parallel to $H_1$ which meets $S\cap H_1'$ in at least 2 points. Assume for example that it is $H_2$. Since $m=3$, these 2 points are included in $L_1$ a line which is a translation of $L$. Consider $H$ the hyperplane containing $L_1$ and $L$. Then, $H$ meets $S\cap H_3$ and $S\cap H_4$ in 1 point (see Figure \ref{fig10a}). So $\#(S\cap H)=7$
\item There are exactly 2 hyperplanes through $L'$ which meets $S$ in 3 points and such that its intersection with $S$ is included in one of the affine hyperplane parallel to $H_1$. 
\\Assume that $H_2$ contains $S\cap \widehat{H}$ where $\widehat{H}$ is the hyperplane through $L'$ different from $H_1$ such that $\#(S\cap \widehat{H})=3$ and $S\cap \widehat{H}$ is included in an hyperplane parallel to $H_1$, say $H_2$. We denote by $L_1=\widehat{H}\cap H_2$. Since for all $H'$ hyperplane $\#(S\cap H')\geq3$, $S\cap H_1'$ meets $H_3$ and $H_4$ in at least one point. Then consider $H$ the hyperplane through $L$ and $L_1$. Since $H$ is different from the hyperplane through $L'$ and $L_1$, $H$ meets $H_3$ and $H_4$ in at least 1 point each (see Figure \ref{fig10b}). So $\#(S\cap H)\geq7$.
\item There are exactly 3 hyperplanes through $L'$ which meets $S$ in 3 points and such that its intersection with $S$ is included in one of the affine hyperplane parallel to $H_1$.
\\If 2 such hyperplanes have their intersection with $S$ included in the same hyperplane parallel to $H_1$, say $H_2$, then $\#(S\cap H_2)\geq7$. Now, assume that they are included in 2 different hyperplanes, $H_2$ and $H_3$. If $S\cap H_1'$ is included in $H_4$ then we consider $H$ the hyperplane through $L$ and $S\cap H_1'$ and $\#(S\cap H)\geq7$. Otherwise, we can assume that $S\cap H_1'$ meets $H_2$ in at least 1 point. Let $H$ be the hyperplane through $L$ and $L_1$ the line containing the minimum weight codeword included in $H_3$. Since $H$ is different from the hyperplane through $L'$ and $L_1$, $H$ meets $S\cap H_2$ in at least 1 point and $\#(S\cap H)\geq7$ (see Figure \ref{fig10c}).
\item  There are 4 hyperplanes through $L'$ which meets $S$ in 3 points and such that its intersection with $S$ is included in one of the affine hyperplane parallel to $H_1$. 
\\If 3 such hyperplanes have their intersection with $S$ included in the same hyperplane parallel to $H_1$, say $H_2$, then $\#(S\cap H_2)\geq7$. Assume that 2 such hyperplanes have their intersection included in the same hyperplane parallel to $H_1$, say $H_2$ and the last one has its intersection with $S$ included in $H_3$. Then, since $\#(S\cap H_4)\geq3$, $\#(S\cap H_1'\cap H_4)\geq3$. \\If $\#(S\cap H_4\cap H_1')=4$, we consider $H$ the hyperplane through $L$ and $S\cap H_1'$ and $\#(S\cap H)\geq7$. Otherwise, there is one point of $S\cap H_4$ included in $H_2$ or $H_3$. If this point is included in $H_2$ then $\#(S\cap H_2)\geq 7$. If it is included in $H_3$, we consider $L_1$ and $L_2$ the 2 lines in $H_2$ containing $S$ which are a translation of $L$. Then either the hyperplane through $L$ and $L_1$ or the hyperplane through $L$ and $L_2$ meets $S\cap H_3$ or $S\cap H_4$ (see Figure \ref{fig10d}). So there is an hyperplane $H$ such that $\#(S\cap H)\geq7$.

Now assume that for each hyperplane $H'$ parallel to $H_1$, there is only one hyperplane through $L'$ which meets $S$ in 3 points such that its intersection with $S$ included in $H'$. If $S\cap H_1'$ is included in an affine hyperplane parallel to $H_1$ then we consider $H$ the hyperplane through $L$ and  $S\cap H_1'$ and $\#(S\cap H)\geq7$. Otherwise, $S\cap H_1'$ meets at least 2 hyperplanes parallel to $H_1$, say $H_2$ and $H_3$ in at least 1 point. For $i\in\{2,3,4\}$, we denote by $H_i'$ the hyperplane through $L'$ such that $S\cap H_i'\subset H_i$. If $\widehat{H}$ the hyperplane through $L$ and $S\cap H'_4$ does not meet $S\cap H_2$ and $S\cap H_3$, then $\widetilde{H}$ the hyperplane through $S\cap H_4'$ and $S\cap H_3'$ meets $S\cap H_2$. Indeed, if $\widehat{H}$ does not meet $S\cap H_2$ we consider 4 hyperplanes through $S\cap H_4'$ different from $H_4$, which intersect $H_2$ in 4 distinct parallel lines. However 2 of these lines meet $S$ (see Figure \ref{fig10e}). So there is an hyperplane $H$ such that $\#(S\cap H)\geq7$.
\end{enumerate} 

\begin{figure}[!h]
\caption{}
\label{fig10}
\begin{center}
\subfloat[Case 1]{\label{fig10a}\begin{tikzpicture}[scale=0.2]
\draw (0,0)--++(5,2)--++(0,8)--++(-5,-2)--cycle;
\draw (6,0)--++(5,2)--++(0,8)--++(-5,-2)--cycle;
\draw (12,0)--++(5,2)--++(0,8)--++(-5,-2)--cycle;
\draw (18,0)--++(5,2)--++(0,8)--++(-5,-2)--cycle;
\draw (5,10) node[above]{$H_1$};
\draw (11,10) node[above]{$H_2$};
\draw (17,10) node[above]{$H_3$};
\draw (23,10) node[above]{$H_4$};
\draw (0,7)--++(5,2);
\draw (1,7+2/5) node{$\bullet$};
\draw (2,7+4/5) node{$\bullet$};
\draw (3,7+6/5) node{$\bullet$};
\draw (0,7) node[left]{$L$};
\draw[dashed] (0,5)--++(5,2);
\draw (0,5) node[left]{$L'$};
\draw (6,5)--++(5,2);
\draw (7,5+2/5) node{$\bullet$};
\draw (9,5+6/5) node{$\bullet$};
\draw (6,5) node[below right]{$L_1$};
\draw (12,1)--++(5,2);
\draw (13,1+2/5) node{$\bullet$};
\draw (12,3)--++(5,2);
\draw (14,3+4/5) node{$\bullet$};
\draw (12,7)--++(5,2);
\draw (16,7+8/5) node{$\bullet$};
\draw (18,1)--++(5,2);
\draw (20,1+4/5) node{$\bullet$};
\draw (18,3)--++(5,2);
\draw (22,3+8/5) node{$\bullet$};
\draw (18,7)--++(5,2);
\draw (21,7+6/5) node{$\bullet$};
\draw[dashed] (12,5)--++(5,2);
\draw[dashed] (18,5)--++(5,2);
\draw (0,5)--(18,5);
\draw (5,7)--(23,7);
\fill[color=gray!60,opacity=0.5] (0,7)--(18,1)--++(5,2)--(5,9)--cycle;
\draw (23,7) node [right]{$H_1'$};
\draw (23,3) node [right]{$H$};
\end{tikzpicture}}\hfill
\subfloat[Case 2]{\label{fig10b}\begin{tikzpicture}[scale=0.2]
\draw (0,0)--++(5,2)--++(0,8)--++(-5,-2)--cycle;
\draw (6,0)--++(5,2)--++(0,8)--++(-5,-2)--cycle;
\draw (12,0)--++(5,2)--++(0,8)--++(-5,-2)--cycle;
\draw (18,0)--++(5,2)--++(0,8)--++(-5,-2)--cycle;
\draw (5,10) node[above]{$H_1$};
\draw (11,10) node[above]{$H_2$};
\draw (17,10) node[above]{$H_3$};
\draw (23,10) node[above]{$H_4$};
\draw (0,7)--++(5,2);
\draw (1,7+2/5) node{$\bullet$};
\draw (2,7+4/5) node{$\bullet$};
\draw (3,7+6/5) node{$\bullet$};
\draw (0,7) node[left]{$L$};
\draw[dashed] (0,5)--++(5,2);
\draw (0,5) node[left]{$L'$};
\draw (6,5)--++(5,2);
\draw (7,5+2/5) node{$\bullet$};
\draw (9,5+6/5) node{$\bullet$};
\draw (10,5+8/5) node{$\bullet$};
\draw (6,5) node[below right]{$L_1$};
\draw (12,1)--++(5,2);
\draw (13,1+2/5) node{$\bullet$};
\draw (12,3)--++(5,2);
\draw (14,3+4/5) node{$\bullet$};
\draw (12,7)--++(5,2);
\draw (16,7+8/5) node{$\bullet$};
\draw (18,1)--++(5,2);
\draw (20,1+4/5) node{$\bullet$};
\draw (18,3)--++(5,2);
\draw (22,3+8/5) node{$\bullet$};
\draw (18,7)--++(5,2);
\draw (21,7+6/5) node{$\bullet$};
\draw[dotted] (12,5)--++(5,2);
\draw[dotted] (18,5)--++(5,2);
\draw (0,5)--(18,5);
\draw (5,7)--(23,7);
\fill[color=gray!60,opacity=0.5] (0,7)--(18,1)--++(5,2)--(5,9)--cycle;
\draw (23,7) node [right]{$\widehat{H}$};
\draw (23,3) node [right]{$H$};
\end{tikzpicture}}\\
\subfloat[Case 3]{\label{fig10c}\begin{tikzpicture}[scale=0.2]
\draw (0,0)--++(5,2)--++(0,8)--++(-5,-2)--cycle;
\draw (6,0)--++(5,2)--++(0,8)--++(-5,-2)--cycle;
\draw (12,0)--++(5,2)--++(0,8)--++(-5,-2)--cycle;
\draw (18,0)--++(5,2)--++(0,8)--++(-5,-2)--cycle;
\draw (5,10) node[above]{$H_1$};
\draw (11,10) node[above]{$H_2$};
\draw (17,10) node[above]{$H_3$};
\draw (23,10) node[above]{$H_4$};
\draw (0,7)--++(5,2);
\draw (1,7+2/5) node{$\bullet$};
\draw (2,7+4/5) node{$\bullet$};
\draw (3,7+6/5) node{$\bullet$};
\draw (0,7) node[left]{$L$};
\draw[dashed] (0,5)--++(5,2);
\draw (0,5) node[left]{$L'$};
\draw[dotted] (6,5)--++(5,2);
\draw (12,5)--++(5,2);
\draw (12,5) node[below right]{$L_1$};
\draw (13,5+2/5) node{$\bullet$};
\draw (15,5+6/5) node{$\bullet$};
\draw (16,5+8/5) node{$\bullet$};
\draw[dotted] (18,5)--++(5,2);
\draw (0,5)--(18,5);
\draw (5,7)--(23,7);
\draw (6,6)--++(5,2);
\draw (10,6+8/5) node{$\bullet$};
\fill[color=gray!60,opacity=0.5] (0,7)--(18,4)--++(5,2)--(5,9)--cycle;
\draw (6,3)--++(5,2);
\draw (10,3+8/5) node{$\bullet$};
\draw (8,3+4/5) node{$\bullet$};
\draw (9,3+6/5) node{$\bullet$};
\draw (6,1)--++(5,2);
\draw (10,1+8/5) node{$\bullet$};
\draw (23,6) node[right]{$H$};
\end{tikzpicture}}\\
\subfloat[Case 4]{\label{fig10d}\begin{tikzpicture}[scale=0.2]
\draw (0,0)--++(5,2)--++(0,8)--++(-5,-2)--cycle;
\draw (6,0)--++(5,2)--++(0,8)--++(-5,-2)--cycle;
\draw (12,0)--++(5,2)--++(0,8)--++(-5,-2)--cycle;
\draw (18,0)--++(5,2)--++(0,8)--++(-5,-2)--cycle;
\draw (5,10) node[above]{$H_1$};
\draw (11,10) node[above]{$H_2$};
\draw (17,10) node[above]{$H_3$};
\draw (23,10) node[above]{$H_4$};
\draw (0,7)--++(5,2);
\draw (1,7+2/5) node{$\bullet$};
\draw (2,7+4/5) node{$\bullet$};
\draw (3,7+6/5) node{$\bullet$};
\draw (0,7) node[left]{$L$};
\draw[dashed] (0,4)--++(5,2);
\draw (0,4) node[left]{$L'$};
\draw (6,3)--++(5,2);
\draw (7,3+2/5) node{$\bullet$};
\draw (10,3+8/5) node{$\bullet$};
\draw (9,3+6/5) node{$\bullet$};
\draw (6,3) node[below right]{$L_2$};
\draw (6,5)--++(5,2);
\draw (8,5+4/5) node{$\bullet$};
\draw (10,5+8/5) node{$\bullet$};
\draw (9,5+6/5) node{$\bullet$};
\draw (11,7) node[above left]{$L_1$};
\draw[dotted] (12,2)--++(5,2);
\draw[dotted] (12,6)--++(5,2);
\draw[dotted] (18,1)--++(5,2);
\draw[dotted] (18,7)--++(5,2);
\draw (0,4)--(18,1);
\draw (0,4)--(18,7);
\draw (5,6)--(23,3);
\draw (5,6)--(23,9);
\draw (12,4)--++(5,2);
\draw (14,4+4/5) node{$\bullet$};
\draw (13,4+2/5) node{$\bullet$};
\draw (16,4+8/5) node{$\bullet$};
\draw (12,3)--++(5,2);
\draw (14,3+4/5) node{$\bullet$};
\draw[dotted] (18,4)--++(5,2);
\draw (0,4)--(18,4);
\draw (5,6)--(23,6);
\draw (18,5)--++(5,2);
\draw (19,5+2/5) node{$\bullet$};
\draw (20,5+4/5) node{$\bullet$};
\draw (21,5+6/5) node{$\bullet$};
\end{tikzpicture}}\hfill
\subfloat[Case 4']{\label{fig10e}\begin{tikzpicture}[scale=0.2]
\draw (0,0)--++(5,2)--++(0,8)--++(-5,-2)--cycle;
\draw (6,0)--++(5,2)--++(0,8)--++(-5,-2)--cycle;
\draw (12,0)--++(5,2)--++(0,8)--++(-5,-2)--cycle;
\draw (18,0)--++(5,2)--++(0,8)--++(-5,-2)--cycle;
\draw (5,10) node[above]{$H_1$};
\draw (11,10) node[above]{$H_2$};
\draw (17,10) node[above]{$H_3$};
\draw (23,10) node[above]{$H_4$};
\draw[dashed] (0,7)--++(5,2);
\draw (0,7) node[left]{$L'$};
\draw (0,1)--++(5,2);
\draw (1,1+2/5) node{$\bullet$};
\draw (2,1+4/5) node{$\bullet$};
\draw (3,1+6/5) node{$\bullet$};
\draw (0,1) node[left]{$L$};
\draw (18,1)--++(5,2);
\draw (22,1+8/5) node{$\bullet$};
\draw (20,1+4/5) node{$\bullet$};
\draw (21,1+6/5) node{$\bullet$};
\draw (23,3) node[right]{$S\cap H_4'$};
\draw (0,1)--(18,1);
\draw (23,3)--(5,3);
\draw( 18,1) node[below left]{$\widehat{H}$};
\draw (12,2)--++(5,2);
\draw (16,2+8/5) node{$\bullet$};
\draw (14,2+4/5) node{$\bullet$};
\draw (13,2+2/5) node{$\bullet$};
\draw (6,3)--++(5,2);
\draw (10,3+8/5) node{$\bullet$};
\draw (8,3+4/5) node{$\bullet$};
\draw (7,3+2/5) node{$\bullet$};
\draw[dotted] (0,4)--+(5,2);
\fill[color=gray!60, opacity=0.5] (0,4)--(18,1)--++(5,2)--(5,6)--cycle;
\draw (0,4) node[ right]{$\widetilde{H}$};
\draw[dotted] (6,5)--++(5,2);
\draw[dotted] (12,3)--++(5,2);
\draw (0,7)--(18,1);
\draw (5,9)--(23,3);
\draw (6,7)--++(5,2);
\draw (9,7+6/5) node{$\bullet$};
\draw (12,7)--++(5,2);
\draw (14,7+4/5) node{$\bullet$};
\end{tikzpicture}}
\end{center}
\end{figure}

In all cases, there exists an affine hyperplane $H$ such that $\#(S\cap H)\geq7$. If $\#(S\cap H)>7$, since $S$ meets all affine hyperplanes in at least 3 points, $\#S>7+3.3=16$ which gives a contradiction. If $\#(S\cap H)=7$, then for all $H'$ parallel to $H$ different form $H$ $\#(S\cap H')=3$. By applying an affine transformation, we can assume that $x_1=0$ is an equation of $H$. Then $g=x_1f\in R_4(3(m-2)+2,m)$ and $|g|=9$. So, $g$ is a second weight codeword of $R_4(3(m-2)+2,m)$ and by Theorem \ref{wpoids2}, the support of $g$ is included in a plane $P$. Since $S$ meets all hyperplanes, $S$ is not included in $P$. Then, $S$ meets all hyperplanes parallel to $P$ in at least 3 points. However $3.3+9=18>16$ which is absurd.
\\\\\indent Now, assume that $m\geq4$. Assume that $S$ is not included in an affine subspace of dimension 3. Then there exists $H$ an affine hyperplane such that $S\cap H$ is not included in a plane and $S$ is not included in $H$. So, by Lemma \ref{cle1}, $S$ meets all affine hyperplanes parallel to $H$ in at least 3 points. 

Assume that for all $H'$ parallel to $H$, $\#(S\cap H')\geq4$, then for reason of cardinality, $\#(S\cap H)=4$. So $S\cap H$ is the support of a second weight codeword of $R_4(3(m-1)+1,m)$ and is included in a plane which is absurd. 
So there exists $H_1$ an affine hyperplane parallel to $H$ such that $\#(S\cap H_1)=3$. Then, $S\cap H_1$ is the support of the minimum weight codeword of $R_4(3(m-1)+1,m)$ and is included in a line $L$. We denote by $\overrightarrow{u}$ a directing vector of $L$ and $a$ the point of $L$ which is not in $S$. 

Let $w_1$, $w_2$, $w_3$ 3 points of $S\cap H$ which are not included in a line. Then, there are at least 2 vectors of $\{\overrightarrow{w_1w_2},\overrightarrow{w_1w_3},\overrightarrow{w_2w_3}\}$ which are not collinear to $\overrightarrow{u}$. Assume that they are $\overrightarrow{w_1w_2}$ and $\overrightarrow{w_1w_3}$. Let $a$ be an affine subspace of codimension 2 included in $H_1$  which contains $a$, $a+\overrightarrow{w_1w_2}$, $a+\overrightarrow{w_1w_3}$ but not $a+\overrightarrow{u}$. Then $S$ does not meet $A$. Assume that $S$ does not meet one hyperplane through $A$. Then $S$ is included in an affine hyperplane parallel to this hyperplane which is absurd by definition of $A$. So, $S$ meets all hyperplanes through $A$ and since $3.4+4=16$, There exists $H_2$ an hyperplane through $A$ such that $\#(S\cap H_2)=4$ and $S\cap H_2$ is included in a plane. For all $H'$ hyperplane through $A$ different from $H_2$, $\#(S\cap H')=3$ and $S\cap H'$ is included in a line. Consider $H'_2$ the hyperplane through $A$ such that $w_1\in H'_2$. Then $w_1$, $w_2$, $w_3\in H'_2$. Since for all $H'$ hyperplane through $A$ different from $H_2$, $S\cap H'$ is included in a line and $w_1$, $w_2$, $w_3$ are not included in a line $H'_2=H_2$. Further more $S\cap H_2$ is included in a plane, so $S\cap H'_2\subset H$.

For all $H'$ hyperplane through $A$ different from $H_2$, $S\cap H'$ is the support of a minimum weight codeword of $R_4(3(m-1)+1,m)$ which does not meet $H_1$, so either $S\cap H'$ is included an affine hyperplane parallel to $H_1$ or $S\cap H'$ meets all affine hyperplanes parallel to $H$ but $H_1$ in 1 point. Since $S\cap H_2$ is included in $H$ and all hyperplanes parallel to $H$ meets $S$ in at least $3$ points, there are only 3 possibilities :
\begin{enumerate}\item For all $H'_2$ hyperplane through $A$, $S\cap H'_2$ is included in an affine hyperplane parallel to $H$.
\item For $H'_2$ hyperplane through $A$ different from $H_2$ and $H_1$, $S\cap H'_2$ meets all affine hyperplanes parallel to $H$ different from $H_1$ in $1$ points.
\item There is $4$ hyperplanes through $A$ such that their intersection with $S$ is included in an affine hyperplane parallel to $H$ and 1 hyperplane through $A$ which meets all hyperplanes parallel to $H$ but $H_1$ in $1$. 
\end{enumerate}
In the two first cases, since $S\cap H$ is not included in a plane and $S$ meets all hyperplanes parallel to $H$ in at least 3 points, $\#(S\cap H)=7$ and for all $H'$ parallel to $H$ different form $H$, $\#(S\cap H')=3$. By applying an affine transformation, we can assume that $x_1=0$ is an equation of $H$. Then $g=x_1f\in R_4(3(m-2)+2,m)$ and $|g|=9$. So, $g$ is a second weight codeword of $R_4(3(m-2)+2,m)$ and by Theorem \ref{wpoids2}, the support of $g$ is included in a plane $P$. Since $S$ is not included in $P$, there exists $H'_1$ an affine hyperplane which contains $P$ but not $S$. Then, $S$ meets all hyperplanes parallel to $H_1'$ in at least 3 points. However $3.3+9=18>16$ which is absurd.

Assume we are in the third case. Since $S\cap H$ is the union of a point and $S\cap H_2$ which is included in a plane and $m\geq4$, there exist $B$ an affine subspace of codimension 2 included in $H$ such that $S$ does not meet $B$ and $S\cap H$ is not included in affine hyperplane parallel to $B$. Then $S$ meets all affine hyperplanes through $B$ in at most 4 points which is a contradiction since $\#(S\cap H)=5$.

So $S$ is included in $G$ an affine subspace of dimension 3. By applying an affine transformation, we can assume that $$G:=\{(x_1,\ldots,x_m)\in \F_q^m:x_4=\ldots=x_m=0\}.$$ Let $g\in B_{3}^q$ defined for all $x=(x_{1},x_2,x_3)\in\F_{q}^{3}$ by $$g(x)=f(x_1,x_2,x_3,0,\ldots,0)$$ and denote by $P\in\F_{q}[X_{1},X_2,X_3]$ its reduced form. Since $$\forall x=(x_1,\ldots,x_m)\in\F_{q}^m, \ f(x)=(1-x_4^{q-1})\ldots(1-x_{m}^{q-1})P(x_{1},x_2,x_3),$$ the reduced form of $f\in R_q(3(m-2)+1,m)$ is $$(1-X_4^{q-1})\ldots(1-X_{m}^{q-1})P(X_{1},X_2,X_3).$$ Then $g\in R_q(4,3)$ and $|g|=|f|=16$. Thus, using the case where $m=3$, we finish the proof of Lemma \ref{m-2}.\end{preuve}

\begin{theoreme}\label{s=1}For $q\geq 4$, $m\geq2$, $0\leq t\leq m-2$, if $f\in R_q(t(q-1)+1,m)$ is such that $|f|=q^{m-t}$, the support of $f$ is an affine subspace of codimension $t$. \end{theoreme}

\begin{preuve}
If $t=0$, the second weight is $q^m$ and we have the result.
\\\\\indent For other cases, we proceed by recursion on $t$. \\\\\indent If $q\geq 5$, we have already proved the case where $t=m-1$ (Theorem \ref{t=m-1}); if $m=2$ and $t=m-2=0$, we have the result. Assume that $m\geq3$.
\\\\\indent For $q=4$, if $m=2$, $t=m-2=0$ and we have the result. If $m\geq3$, we have already proved the case $t=m-2$ (Lemma \ref{m-2}). Furthermore, if $m=3$ we have considered all cases. Assume $m\geq4$
\\\\\indent Let $1\leq t\leq m-2$ (or for $q=4$, $1\leq t\leq m-3$). Assume that the support of $f\in R_q((t+1)(q-1)+1,m)$ such that $|f|=q^{m-t-1}$ is an affine subspace of codimension $t+1$.
\\\\\indent Let $f\in R_q(t(q-1)+1,m)$ such that $|f|=q^{m-t}$. We denote by $S$ its support. Assume that $S$ is not included in an affine subspace of codimension $t$. Then there exists $H$ an affine hyperplane such that $S  \cap H$ is not included in an affine subspace of codimension $t+1$ and $S\cap H\neq S$. Then, by Lemma \ref{cle1}, $S$ meets all affine hyperplanes parallel to $H$ and for all $H'$ hyperplane parallel to $H$, $$\#(S\cap H')\geq (q-1)q^{m-t-2}.$$ 
If for all $H'$ hyperplane parallel to $H$, $\#(S\cap H')>(q-1)q^{m-t-2}$ then, for reason of cardinality, $\#(S\cap H)=q^{m-t-1}$. So $S\cap H$ is the support of a second weight codeword of $R_q((t+1)(q-1)+1,m)$ and is included in an affine subspace of codimension $t+1$ which is a contradiction.
\\\indent So there exists $H_1$ parallel to $H$ such that $\#(S\cap H_1)=(q-1)q^{m-t-2}$. Then $S\cap H_1$ is the support of a minimal weight codeword of $R_q((t+1)(q-1)+1,m)$. Hence, $S\cap H_1$ is the union of $q-1$ affine subspaces of codimension $t+2$ included in an affine subspace of codimension $t+1$. 
\\\indent Let $A$ be an affine subspace of codimension $2$ included in $H_1$ such that $A$ meets the affine subspace of codimension $t+1$ which contains $S\cap H_1$ in the affine subspace of codimension $t+2$ which does not meet $S$. Assume that there is an affine hyperplane through $A$ which does not meet $S$. Then, by Lemma \ref{cle1}, $S$ is included in an affine hyperplane parallel to this hyperplane which is absurd by construction of $A$. So, $S$ meets all hyperplanes through $A$. Furthermore, $$q^{m-t}=q^{m-t-1}+q(q-1)q^{m-t-2}.$$
So $S$ meets one of the hyperplane through $A$ in $q^{m-t-1}$ points, say $H_2$, and all the others in $(q-1)q^{m-t-2}$ points. \\Since $H_2\neq H_1$, $H_2\cap H_1=A$ and $S\cap H_2\cap H_1=\emptyset$. So, $S\cap H_2$ is the support of a second weight codewords of $R_q((t+1)(q-1)+1,m)$ which does not meet $H_1$. Hence, $S\cap H_2$ is included in one of the affine hyperplanes parallel to $H$. Furthermore, for all $H_2'$ hyperplane through $A$ different from $H_2$ and $H_1$, $S\cap H_2'$ is the support of a minimum weight codeword of $R_q((t+1)(q-1)+1,m)$ which does not meet $H_1$, so it meets all affine hyperplanes parallel to $H_1$ different from $H_1$ in $q^{m-t-2}$ points or is included in an affine hyperplane parallel to $H_1$. Since $S\cap H_2$ is included in one of the affine hyperplanes parallel to $H$ and all hyperplanes parallel to $H$ meets $S$ in at least $(q-1)q^{m-t-2}$ points, there are only 3 possibilities :
\begin{enumerate}\item For all $H'_2$ hyperplane through $A$, $S\cap H'_2$ is included in an affine hyperplane parallel to $H$.
\item For $H'_2$ hyperplane through $A$ different from $H_2$ and $H_1$, $S\cap H'_2$ meets all affine hyperplanes parallel to $H$ different from $H_1$ in $q^{m-t-2}$ points.
\item There is $q$ hyperplanes through $A$ such that their intersection with $S$ is included in an affine hyperplane parallel to $H$ and 1 hyperplane through $A$ which meets all hyperplanes parallel to $H$ but $H_1$ in $q^{m-t-2}$. 
\end{enumerate}
In the two first cases, if $S\cap H_2$ is not included in $H'$ parallel to $H$, \\$\#(S\cap H')=(q-1)q^{m-t-2}$ and $S\cap H'$ is the support of a minimum weight codewords of $R_q((t+1)(q-1)+1,m)$. So $S\cap H'$ is included in an affine subspace of codimension $t+1$. Then, necessarily, $S\cap H_2$ is included in $H$. For all $H'$ parallel to $H$ but $H$, $\#(S\cap H')=(q-1)q^{m-t-2}$. In the third case, for all $H'$ hyperplane parallel to $H$ different from $H_1$ which does not contain $S\cap H_2$, $\#(S\cap H')=q^{m-t-1}$. So $S\cap H'$ is the support of a second weight codeword of $R_q((t+1)(q-1)+1,m)$ and is an affine subspace of dimension $m-t-1$. Then, $S\cap H_2\subset H$ and $\#(S\cap H)=q^{m-t-2}+q^{m-t-1}$, $\#(S\cap H_1)=(q-1)q^{m-t-2}$. So if we are in the last case for reason of cardinality, for all $A'$ affine subspace of codimension 2 included in $H_1$ such that $A'$ meets the affine subspace of codimension $t+1$ which contains $S\cap H_1$ in the affine subspace of codimension $t+2$ which does not meet $S$ we are also in case 3. Then $S$ is the union of affine subspaces of dimension $m-t-2$ which are a translation of the affine subspace of codimension $t+2$ which does not meet $S$ in $S\cap H_1$. Then, since $S\cap H_2$ is the support of a second weight codeword of $R_q((t+1)(q-1)+1,m)$, it is an affine subspace of dimension $m-t-1$. So $S\cap H$ is the union of an affine subspace of dimension $m-t-1$ and an affine subspace of dimension $m-t-2$. Since $S$ is the union of affine subspaces of dimension $m-t-2$ which are a translation of an affine subspace of codimension $t+2$, there exists $B$ an affine subspace of codimension 2 such that $B$ does not meet $S$ and $S\cap H$ is not included in an affine subspace of codimension 2 parallel to $B$. Now, we consider all affine hyperplanes through $B$. Assume that there exists $G$ an affine hyperplane through $B$ which does not meet $S$. Then, $S$ is included in an affine hyperplane parallel to $G$ which is absurd by construction of $B$. So, $S$ meets all hyperplanes through $B$ and there exists $G_1$ hyperplane through $B$ such that $\#(S\cap G_1)=q^{m-t-1}$ and for all $G$ through $B$ but $G_1$, $\#(S\cap G)=(q-1)q^{m-t-2}$ which is absurd since $\#(S\cap H)=q^{m-t-1}+q^{m-t-2}$.
Finally, we are in case 1 or 2.   
\\\indent By applying an affine transformation, we can assume that $x_1=0$ is an equation of $H$. Now, consider $g$ the function defined by $$\forall x=(x_1,\ldots,x_m)\in\F_q^m\quad g(x)=x_1f(x).$$
Then $\deg(g)\leq t(q-1)+2$ and $|g|=(q-1)^2q^{m-t-2}$. So, $g$ is a second weight codeword of $R_q(t(q-1)+2,m)$ and by Theorem \ref{wpoids2}, the support of $g$ is included in an affine subspace of codimension $t$. 
\\\indent Let $H_3$ be an affine hyperplane containing the support of $g$ but not $S$. Then, $\#(S\cap H_3)\geq (q-1)^2q^{m-t-2}$.
Furthermore, since $S\not\subset H_3$, $S$ meets all affine hyperplanes parallel to $H_3$ in at least $(q-1)q^{m-t-2}$. Finally, $$\#S\geq2(q-1)^2q^{m-t-2}>q^{m-t}.$$ 
\\We get a contradiction. So $S$ is included in an affine subspace of codimension $t$. For reason of cardinality, $S$ is an affine subspace of codimension $t$. \end{preuve}

\subsection{Case where $q=3$, proof of Theorem \ref{1q=3}}

\begin{lemme}\label{31}Let $m\geq2$, $0\leq t\leq m-2$, $f\in R_3(2t+1,m)$ such that \\$|f|=8.3^{m-t-2}$. If $H$ is an affine hyperplane of $\F_q^m$ such that $S\cap H\neq \emptyset$ and $S\cap H\neq S$ then either $S$ meets 2 hyperplanes parallel to $H$ in $4.3^{m-t-2}$points or $S$ meets all affine hyperplanes parallel to $H$.\end{lemme}

\begin{preuve}By applying an affine transformation, we can assume that $x_1=0$ is an equation of $H$. We denote by $H_a$ the affine hyperplanes parallel to $H$ of equation $x_1=a$, $a\in\F_q$. Let $I:=\{a\in\F_q : S\cap H_a=\emptyset\}$ and $k:=\#I$. Since $S\cap H\neq \emptyset$ and $S\cap H\neq S$, $k\leq q-2=1$. Assume $k=1$. For all $c\not\in I$ we define 
$$\forall x=(x_1,\ldots,x_m)\in\F_q^m,\quad f_c(x)=f(x)\prod_{a\not\in I,a\neq c}(x_1-a).$$
Then $\deg(f_c)=(t+1)2$ and $|f_c|\geq 3^{m-t-1}$. Assume that there exists $H'$ an affine hyperplane parallel to $H$ such that $\#(S\cap H')=3^{m-t-1}$ and $S\cap H'$ is the support of a minimal weight codeword of $R_3(2(t+1),m)$. Then consider $A$ an affine subspace of codimension 2 included in $H'$ containing $S\cap H'$ and $A'$ an affine subspace of codimension 2 included in $H'$ parallel to $A$. We denote by $k$ the number of hyperplanes through $A$ which meet $S$ and by $k'$ the number of affine hyperplanes through $A'$ which meet $S$ in $3^{m-t-1}$ points. Then $$k'3^{m-t-1}+(k-k')4.3^{m-t-2}\leq 8.3^{m-t-2}.$$ Since $\#S>\#(S\cap H')$ and $k'\leq k$, we get $k=2$. Then, if we denote by $H''$ the other hyperplane parallel to $H'$ which meets $S$, $S\cap H''$ is included in an affine subspace of codimension 2 which is a translation of $A$. By applying this argument to all affine subspaces of codimension 2 included in $H'$ and containing $S\cap H'$, we get that $S\cap H''$ is included in a an affine subspace of dimension $m-t-1$. For reason of cardinality this is absurd. If $|f_c|>3^{m-t-1}$ then $|f_c|\geq 4.3^{m-t-2}$. For reason of cardinality, we have the result.\end{preuve}

Now, we prove Proposition \ref{1q=3}.
\begin{itemize}\item 
First, we prove the case where $t=1$. Obviously, $S$ is included in an affine subspace of dimension $m$. Assume that $S$ meets all affine hyperplanes of $\F_q^m$. Then for all $H'$ affine hyperplane of $\F_q^m$, $\#(S\cap H')\geq 2.3^{m-3}$ and by Lemma \ref{inter}, there exists $H$ an affine hyperplane such that $$\#(S\cap H)=2.3^{m-3}.$$ Then $S\cap H$ is the support of a minimum weight codeword of $R_3(5,m)$. So it is the union of $P_1$, $P_2$ 2 parallel affine subspaces of dimension $m-3$ included in an affine subspace of dimension $m-2$. Let $A$ be an affine subspace of codimension 2 included in $H$, containing $P_1$ and different from the affine subspace of codimension 2 containing $S\cap H$. Then there exists $A'$ an affine hyperplane of codimension 2 included in $H$ parallel to $A$ which does not meet $S$. We denote by $k$ the number of affine hyperplanes through $A'$ which meet $S$ in $2.3^{m-3}$ points. Then, if $m\geq4$, $$k2.3^{m-3}+(4-k)8.3^{m-4}\leq 8.3^{m-3}$$ which means that $k\geq 4$. If $m=3$, $2k+(4-k)3\leq 8$ which also means that $k\geq 4$. Then for all $H'$ hyperplane through $A$ different from $H$, $S\cap H'$ is a minimal weight codeword of $R_3(5,m)$ which does not meet $H$ and either $S\cap H'$ is included in one of the hyperplanes parallel to $H$ or $S\cap H'$ meets the 2 hyperplanes parallel to $H$ different from $H$. In all cases, $S$ is the union of 8 affine subspace of dimension $m-3$. By applying this argument to all affine subspaces of codimension 2 included in $H$, containing $P_1$ and different from the affine subspace of codimension 2 containing $S\cap H$, we get that these 8 affine subspaces are a translation of $P_1$.

Choose $H_1$ one of the hyperplanes through $A'$ and consider $H_2$ and $H_3$ the 2 hyperplanes parallel to $H_1$. Since $\#(S\cap H_1)=2.3^{m-3}$ and $S$ meets all hyperplanes in at least $2.3^{m-3}$ points, either $\#(S\cap H_2)=3.3^{m-3}$ and $\#(S\cap H_3)=3.3^{m-3}$ or $\#(S\cap H_2)=2.3^{m-3}$ and $\#(S\cap H_3)=4.3^{m-3}$. 

First consider the case where $\#(S\cap H_2)=3.3^{m-3}$ and \\$\#(S\cap H_3)=3.3^{m-3}$. Then there exists an affine subspace of codimension 2 in $H_2$ which does not meet $S$. We denote by $k'$ the number of hyperplanes through $A$ which meet $S$ in $2.3^{m-3}$ points. Then , we have $k'\geq4$ which is absurd since $\#(S\cap H_2)=3.3^{m-3}$.

Now, consider the case where $\#(S\cap H_2)=2.3^{m-3}$ and \\$\#(S\cap H_3)=4.3^{m-3}$. By applying an affine transformation, we can assume that $x_1=0$ is an equation of $H_3$. Then $x_1.f$ is a codeword of $R_3(4,m)$ and $|x_1.f|=4.3^{m-3}$. So, by Theorem \ref{0}, its support is included in an affine hyperplane $H_1'$ and $S\cap H_1'\cap H_3=\emptyset$. So $S$ is included $H_1'$ and $H_3$ and there exists an affine hyperplane through $H_1'\cap H_3$ which does not meet $S$ which is absurd.

Finally there exists an affine hyperplane $G_1$ which does not meet $S$. So, by Lemma \ref{31}, $S$ meets $G_2$ and $G_3$ the 2 hyperplanes parallel to $G_1$ in $4.3^{m-3}$ points. Then, Theorem \ref{0}, $G_2\setminus S$ and $G_3\setminus S$ are the union of two non parallel affine subspaces of codimension 2. Consider $A$ one of the affine subspaces of codimension 2 in $G_2\setminus S$. Assume that all hyperplanes through $A$ meet $S$. So for all $G'$ hyperplane through $A$, $\#(G'\setminus S)\leq 7.3^{m-3}$.  Furthermore, one of the hyperplanes through $A$, say $G$, meets $G_3\setminus S$ in at least $2.3^{m-3}$, then $\#(G\setminus S)\geq 2.3^{m-2}+2.3^{m-3}$ which is absurd (see Figure \ref{fig11}). So there exists $G'$ through $A$ which does not meet $S$. By applying the same argument to the other affine subspace of dimension 2 of $G_2\setminus S$, we get the result for $t=1$.
\begin{figure}[!h]
\caption{}
\label{fig11}
\begin{center}\begin{tikzpicture}[scale=0.2]
\draw (0,0)--++(4,2)--++(0,8)--++(-4,-2)--cycle;
\draw[dotted] (5,0)--++(4,2)--++(0,8)--++(-4,-2)--cycle;
\draw (10,0)--++(4,2)--++(0,8)--++(-4,-2)--cycle;
\draw (4,10) node[above]{$G_2$};
\draw (9,10) node[above]{$G_1$};
\draw (14,10) node[above]{$G_3$};
\draw (0,2)--++(4,2);
\draw (0,6)--++(4,2);
\draw[dotted] (0,4)--++(4,2);
\draw (10,4)--++(4,-1);
\draw[dotted] (10,11/2)--++(4,-1);
\draw (10,7)--++(4,-1);
\draw[dotted](2,1)--++(0,8);
\draw[white] (2,3) node{$\bullet$};
\draw[white] (2,7) node{$\bullet$};
\draw[dotted] (10,0)--(14,10);
\draw[white] (138/11,{-138/11+38)/4}) node{$\bullet$};
\draw[white] (126/11,{-126/11+26)/4}) node{$\bullet$};
\end{tikzpicture}\end{center}
\end{figure}

\item We prove by recursion on $t$ that $S$ is included in an affine subspace of dimension $m-t+1$. Consider first the case where $t=m-2$. If $m=3$ then $t=1$ and we have already consider this case. Assume that $m\geq4$. Let $f\in R_3(2(m-2)+1,m)$ such that $|f|=8$. Assume that $S$ is not included in an affine subspace of dimension 3. Let $w_1$, $w_2$, $w_3$, $w_4$ 4 points of $S$ which are not included in a plane. Since $S$ is not included in an affine subspace of dimension 3, there exists $H$ an affine hyperplane such that $H$ contains $w_1$, $w_2$, $w_3$, $w_4$ but $S$ is not included in $H$. Then by Lemma \ref{31} either $S$ meets 2 hyperplanes parallel to $H$ in 4 points or $S$ meets all hyperplanes parallel to $H$. 

If $S$ meets 2 hyperplanes parallel to $H$ then $S\cap H$ is the support of a second weight codeword of $R_3(2(m-1),m)$ so is included in a plane which is absurd since $w_1$, $w_2$, $w_3$, $w_4\in S\cap H$. So $S$ meets all hyperplanes parallel to $H$ and for all $H'$ hyperplane parallel to $H$, $\#(S\cap H')\geq 2$. Since $\#S=8$ and $\#(S\cap H)\geq 4$, for all $H'$ hyperplane parallel to $H$ different from $H$ $\#(S\cap H')=2$ and $\#(S\cap H)=4$. 

By applying an affine transformation, we can assume that $x_1=0$ is an equation of $H$. Then $x_1.f\in R_3(2(m-1),m)$ and $|x_1.f|=4$ so $x_1.f$ is a second weight codeword of $R_3(2(m-1),m)$ and its support is included in a plane $P$ not included in $H$. Let $H'$ be an affine hyperplane which contains $P$ and a point of $(S\cap H)\setminus P$ but not all the points of $S\cap H$. Then, $\#(S\cap H')\geq5$ and $S\cap H'\neq S$. By applying the same argument to $H'$ than to $H$ we get a contradiction for reason of cardinality. 
\item If $m\leq 4$, we have already considered all the possible values for $t$. Assume that $m\geq5$. Let $2\leq t\leq m-3$. Assume that if $f\in R_3(2(t+1)+1,m)$ is such that $|f|=8.3^{m-t-3}$ then its support is included in an affine subspace of dimension $m-t$. Let $f\in R_3(2t+1,m)$ such that $|f|=8.3^{m-t-2}$ and denote by $S$ its support. Assume that $S$ is not included in an affine subspace of dimension $m-t+1$. Then, there exists $H$ an affine hyperplane such that $S\cap H\neq S$ and $S\cap H$ is not included in an affine subspace of dimension $m-t$. So, by Lemma \ref{31}, either $S$ meets 2 affine hyperplanes parallel to $H$ in $4. 3^{m-t-2}$ points or $S$ meets all affine hyperplanes parallel to $H$. 

If $S$ meets 2 affine hyperplanes in $4.3^{m-t-2}$ points, $S\cap H$ is the support of a second weight codeword of $R_3(2(t+1),m)$ and is included in an affine subspace of dimension $m-t$ which is absurd. So $S$ meets all affine hyperplanes parallel to $H$ and for all $H'$ hyperplane parallel to $H$, $$\#(S\cap H')\geq 2.3^{m-t-2}.$$ Assume that for all $H'$ parallel to $H$, $\#(S\cap H')>2.3^{m-t-2}$. Then, for reason of cardinality $\#(S\cap H)=8.3^{m-t-3}$ and $S\cap H$ is the support of a second weight codeword of $R_3(2(t+1)+1,m)$ which is absurd since $S\cap H$ is not included in an affine subspace of dimension $m-t$. So there exists $H_1$ parallel to $H$ such that $\#(S\cap H_1)=2.3^{m-t-2}$ and $S\cap H_1$ is the support of a minimal weight codeword of $R_3(2(t+1)+1,m)$ so $S\cap H_1$ is the union of $P_1$ and $P_2$ 2 parallel affine subspaces of dimension $m-t-2$ included in an affine subspace of dimension $m-t-1$. 

Let $A$ be an affine subspace of codimension 2 included in $H_1$ and containing $P_1$ and such that $A\cap P_2=\emptyset$. Let $A'$ be an affine subspace of codimension 2 included in $H_1$ parallel to $A$ which does not meet $S$. Assume that there exists $H'_1$ an affine hyperplane through $A'$ which does not meet $S$. Then, $S$ meets $H_2'$ and $H_3'$ the 2 hyperplanes parallel to $H_1'$ different from $H_1'$ in $4.3^{m-t-2}$ points. For example, we can assume that $A\subset H_2'$. Then, $S\cap H_3'$ is the support of a second weight codeword of $R_3(2(t+1),m)$. So $S\cap H_3'$ meets $H$ in 0, $3^{m-t-2}$, $2.3^{m-t-2}$ or $4.3^{m-t-2}$ points. Since $S$ meets all hyperplanes parallel to $H$ in at least $2.3^{m-t-2}$ points, if $$\#(S\cap H\cap H_3')=4.3^{m-t-2},$$ $S\cap H\cap H_2'=\emptyset$. So $S\cap H$ is included in an affine subspace of dimension $m-t$ which is absurd. So $S\cap H_2'$ and $S\cap H_3'$ are the support of second weight codewords of $R_3(2(t+1),m)$ not included in $H$, then their intersection with $H$ is the union of at most 2 disjoint affine subspaces of dimension $m-t-2$.

Now assume that $S$ meets all hyperplanes through $A'$. We denote by $k$ the number of the hyperplanes through $A$ which meet $S$ in $2.3^{m-t-2}$ points. Then
$$k2.3^{m-t-2}+(4-k)8.3^{m-t-3}\leq 8.3^{m-t-2}$$ which means that $k\geq4$. So for all $H'$ affine hyperplane through $A'$ different from $H_1$, $S\cap H'$ is the support of minimum weight codeword of $R_3(2(t+1)+1,m)$ which does not meet $H_1$. So either $S\cap H'$ is included in $H$ or $S\cap H'$ meets $S$ in an affine subspace of dimension $m-t-2$. In both cases , $S\cap H$ is the union of at most 4 disjoint affine subspaces of dimension $m-t-2$.
By applying this argument to all affine subspaces of dimension 2 included in $H_1$ containing $P_1$ but not $P_2$, we get that $S\cap H$ is the union of 4 affine subspaces of dimension $m-t-2$ which are a translation of $P_1$. This gives a contradiction since $S\cap H$ is not included in an affine subspace of dimension $m-t$. So $S$ is included in an affine subspace of dimension $m-t+1$. 
\item Let $f\in R_3(2t+1,m)$ such that $|f|=8.3^{m-t-2}$ and $A$ the affine subspace of dimension $m-t+1$ containing $S$. By applying an affine transformation, we can assume $$A:=\{(x_1,\ldots,x_m)\in\F_q^m:x_1=\ldots=x_{t-1}=0\}.$$ 
Let $g\in B_{m-t+1}^3$ defined for all $x=(x_{t},\ldots,x_m)\in\F_{3}^{m-t+1}$ by $$g(x)=f(0,\ldots,0,x_{t},\ldots,x_m)$$ and denote by $P\in\F_{3}[X_{t},\ldots,X_m]$ its reduced form. Since $$\forall x=(x_1,\ldots,x_m)\in\F_{3}^m, \ f(x)=(1-x_1^{2})\ldots(1-x_{t-1}^{2})P(x_{t},\ldots,x_m),$$ the reduced form of $f\in R_3(t(q-1)+s,m)$ is $$(1-X_1^{2})\ldots(1-X_{t-1}^{2})P(X_{t},\ldots,X_m).$$ Then $g\in R_3(3,m-t+1)$ and $|g|=|f|=8.3^{m-t-2}$. Thus, using the case where $t=1$, we finish the proof of Proposition \ref{1q=3}.
\end{itemize}

\appendix
\section{Appendix : Blocking sets}

Blocking sets have been studied by Bruen in \cite{MR0303406,MR0251629,MR2766082} in the case of projective planes. Erickson extends his results to affine planes in \cite{erickson1974counting}.

\begin{defi}Let $S$ be a subset of the affine space $\F_q^2$. We say that $S$ is a blocking set of order $n$ of $\F_q^2$ if for all line $L$ in $\F_q^2$, $\#(S\cap L)\geq n$ and $\#((\F_q^2\setminus S)\cap L)\geq n$.\end{defi}

\begin{proposition}[Lemma 4.2 in \cite{erickson1974counting}]\label{4.2}Let $q\geq 3$, $1\leq b\leq q-1$ and \\$f\in R_q(b,2)$. If  $f$ has no linear factor and $|f|\leq (q-b+1)(q-1)$, then the support of $f$ is a blocking set of order $(q-b)$ of $\F_q^2$.\end{proposition}

In \cite{erickson1974counting} Erickson make the following conjecture. It has been proved by Bruen in \cite{MR2766082}.

\begin{theoreme}[Conjecture 4.14 in \cite{erickson1974counting}]\label{4.14}If $S$ is a blocking set of order $n$ in $\F_q^2$, then $\#S\geq nq+q-n$.\end{theoreme}
\bibliographystyle{plain}

\bibliography{C:/Users/Elodie/Dropbox/bibliothese}

\end{document}